\documentclass{amsart}
\usepackage{graphicx}
\usepackage{amsmath}
\usepackage{amsthm}
\usepackage{amsfonts}
\usepackage{amssymb, amscd}

\usepackage{color}
\usepackage[all]{xy}
\newtheorem{thm}{Theorem}[section]
\newtheorem{cor}[thm]{Corollary}
\newtheorem{lem}[thm]{Lemma}
\newtheorem{proposition}[thm]{Proposition}
\newtheorem{example}[thm]{Example}
\theoremstyle{definition}
\newtheorem{definition}[thm]{Definition}
\theoremstyle{remark}
\newtheorem{remark}[thm]{Remark}
\numberwithin{equation}{section}
\newcommand{\K}{\mathbb K}

\newcommand{\dl}{\displaystyle}

\begin{document}

\title[  Quadratic $n$-ary Hom-Nambu algebras ]
{Quadratic $n$-ary Hom-Nambu algebras}%
\author{F. AMMAR, S. MABROUK and A. MAKHLOUF  }%
\address{Abdenacer Makhlouf, Universit\'{e} de Haute Alsace,  Laboratoire de Math\'{e}matiques, Informatique et Applications,
4, rue des Fr\`{e}res Lumi\`{e}re F-68093 Mulhouse, France}%
\email{Abdenacer.Makhlouf@uha.fr}
\address{Faouzi Ammar and  Sami  Mabrouk, Universit\'{e} de Sfax,  Facult\'{e} des Sciences, Sfax Tunisia}%
\email{Faouzi.Ammar@fss.rnu.tn, \ Mabrouksami00@yahoo.fr}

%\thanks {}

\subjclass[2000]{17A30,17A36,17A42,17A45}
\keywords{$n$-ary Nambu algebra, $n$-ary Nambu-Lie algebra, $n$-ary Hom-Nambu algebra, $n$-ary Hom-Nambu-Lie algebra, Quadratic $n$-ary Hom-Nambu algebra, Representation, centroid}
\date{}
%
%\dedicatory{}%
%\commby{}%
% ----------------------------------------------------------------
\begin{abstract}
The purpose  of this paper is to introduce and study quadratic $n$-ary Hom-Nambu algebras, which are  $n$-ary Hom-Nambu algebras with an invariant, nondegenerate and symmetric bilinear forms that are also $\alpha$-symmetric and $\beta$-invariant where $\alpha$ and $\beta$ are twisting maps. We provide constructions of these $n$-ary algebras by using twisting principles, tensor product and T*-extension. Also is discussed their connections with representation theory and centroids. Moreover we show that one may derive from quadratic $n$-ary Hom-Nambu algebra ones of increasingly higher arities and that under suitable assumptions it reduces to a quadratic  $(n-1)$-ary  Hom-Nambu algebra.
\end{abstract}
\maketitle
% ----------------------------------------------------------------

\section*{Introduction}
The main motivations to study $n$-ary algebras came firstly from Nambu mechanics \cite{Nam} where a ternary bracket allows to use  more than one hamiltonian and recently from string theory and M-branes which involve
naturally an algebra with ternary operation called Bagger-Lambert algebra \cite{bag}. Also ternary operations appeared in the study of some quarks models see  \cite{ker1,ker2, ker3}.

Algebras endowed with invariant nondegenerate symmetric bilinear form (scalar product) appeared also naturally in several domains in mathematics and physics. Such algebras were intensively studied for binary Lie and associative algebras. The main result is that called double extension given by Medina and Revoy \cite{medina} or T*-extension given by Bordemann \cite{Bord}. These fundamental results were extended to $n$-ary algebra in \cite{metricalgebras}. The extension to Hom-setting for binary case was introduced and studied in \cite{BenayadiMakhlouf}.

In this paper we summarize in the first Section  definitions of $n$-ary Hom-Nambu  algebras and recall the constructions  using twisting principles and tensor product with $n$-ary algebras of Hom-associative type.
In Section $2$, we introduce the notion of quadratic $n$-ary Hom-Nambu algebra, generalizing the notion introduced for binary Hom-Lie algebras in \cite{BenayadiMakhlouf}. A more general notion called Hom-quadratic $n$-ary Hom-Nambu algebra is introduced by twisting the invariance identity. In Section 3, we show that a quadratic $n$-ary Hom-Nambu algebra gives rise a  quadratic Hom-Leibniz algebra. A connection with representation theory is discussed in Section 4. We deal in particular with adjoint and coadjoint representations, extending  the representation  theory  initiated in \cite{BenayadiMakhlouf,Sheng}. Several procedures to built quadratic $n$-ary Hom-Nambu algebras are provided in Section 5. We use twisting principles, tensor product and T*-extension to construct quadratic $n$-ary Hom-Nambu algebras. Moreover we show that one may derive from quadratic $n$-ary Hom-Nambu algebra ones of increasingly higher arities and that under suitable assumptions it reduces to a quadratic  $(n-1)$-ary  Hom-Nambu algebra. Also real  Faulkner construction is used to obtain ternary Hom-Nambu algebras. The last Section is dedicated to introduce and study the centroids of $n$-ary Hom-Nambu algebras and their properties. We supply a construction procedure of quadratic $n$-ary Hom-Nambu algebras using elements of the centroid.

Most of the results concern $n$-ary Hom-Nambu algebras. Naturally they are valid and may be stated for $n$-ary Hom-Nambu-Lie algebras.
\section{The n-ary Hom-Nambu algebras}
  Throughout this paper, we will for simplicity of exposition assume that $\mathbb{K}$ is an algebraically closed
field of characteristic zero, even though for most of the general definitions and results in the paper this
assumption is not essential. 
\subsection{Definitions}
In this section, we recall  definitions of $n$-ary Hom-Nambu algebras and $n$-ary Hom-Nambu-Lie algebras, generalizing $n$-ary Nambu algebras and $n$-ary Nambu-Lie algebras
(called also Filippov algebras). They were  introduced in \cite{AMS} by Ataguema, Makhlouf and Silvestrov.
\begin{definition}
An \emph{$n$-ary Hom-Nambu} algebra is a triple $(\mathcal{N}, [\cdot ,..., \cdot],  \widetilde{\alpha} )$ consisting of a vector space  $\mathcal{N}$, an
$n$-linear map $[\cdot ,..., \cdot ] :  \mathcal{N}^{ n}\longrightarrow \mathcal{N}$ and a family
$\widetilde{\alpha}=(\alpha_i)_{1\leq i\leq n-1}$ of  linear maps $ \alpha_i:\ \ \mathcal{N}\longrightarrow \mathcal{N}$, satisfying \\
  \begin{eqnarray}\label{NambuIdentity}
  && \big[\alpha_1(x_1),....,\alpha_{n-1}(x_{n-1}),[y_1,....,y_{n}]\big]= \\ \nonumber
&& \sum_{i=1}^{n}\big[\alpha_1(y_1),....,\alpha_{i-1}(y_{i-1}),[x_1,....,x_{n-1},y_i]
  ,\alpha_i(y_{i+1}),...,\alpha_{n-1}(y_n)\big],
  \end{eqnarray}
  for all $(x_1,..., x_{n-1})\in \mathcal{N}^{ n-1}$, $(y_1,...,  y_n)\in \mathcal{N}^{ n}.$\\
  The identity \eqref{NambuIdentity} is called \emph{Hom-Nambu identity}.
  \end{definition}

Let
$x=(x_1,\ldots,x_{n-1})\in \mathcal{N}^{n-1}$, $\widetilde{\alpha}
(x)=(\alpha_1(x_1),\ldots,\alpha_{n-1}(x_{n-1}))\in \mathcal{N}^{n-1}$ and
$y\in \mathcal{N}$. We define an adjoint map  $ad(X)$ as  a linear map on $\mathcal{N}$,
such that
\begin{equation}\label{adjointMapNaire}
ad(x)(y)=[x_{1},\cdots,x_{n-1},y].
\end{equation}

Then the Hom-Nambu identity \eqref{NambuIdentity} may be written in terms of adjoint map as
\begin{equation*}
ad(\widetilde{\alpha} (x))( [x_{n},...,x_{2n-1}])=
\sum_{i=n}^{2n-1}{[\alpha_1(x_{n}),...,\alpha_{i-n}(x_{i-1}),
ad(x)(x_{i}), \alpha_{i-n+1}(x_{i+1}) ...,\alpha_{n-1}(x_{2n-1})].}
\end{equation*}

\begin{remark}
When the maps $(\alpha_i)_{1\leq i\leq n-1}$ are all identity maps, one recovers the classical $n$-ary Nambu algebras. The Hom-Nambu Identity \eqref{NambuIdentity}, for $n=2$,  corresponds to Hom-Jacobi identity (see \cite{MS}), which reduces to Jacobi identity when $\alpha_1=id$.
\end{remark}

Let $(\mathcal{N},[\cdot,\dots,\cdot],\widetilde{\alpha})$ and
$(\mathcal{N}',[\cdot,\dots,\cdot]',\widetilde{\alpha}')$ be two $n$-ary Hom-Nambu
algebras  where $\widetilde{\alpha}=(\alpha_{i})_{i=1,\cdots,n-1}$ and
$\widetilde{\alpha}'=(\alpha'_{i})_{i=1,\cdots,n-1}$. A linear map $f:
\mathcal{N}\rightarrow \mathcal{N}'$ is an  $n$-ary Hom-Nambu algebras \emph{morphism}  if it satisfies
\begin{eqnarray*}f ([x_{1},\cdots,x_{n}])&=&
[f (x_{1}),\cdots,f (x_{n})]'\\
f \circ \alpha_i&=&\alpha'_i\circ f \quad \forall i=1,n-1.
\end{eqnarray*}

\begin{definition}
An $n$-ary Hom-Nambu algebra $(\mathcal{N}, [\cdot ,..., \cdot],  \widetilde{ \alpha} )$ where  $\widetilde{\alpha}=(\alpha_i)_{1\leq i\leq n-1}$
is called \emph{$n$-ary Hom-Nambu-Lie algebras} algebra if the bracket is skew-symmetric that is
\begin{equation}
[x_{\sigma(1)},..,x_{\sigma(n)}]=Sgn(\sigma)[x_1,..,x_n],\ \ \forall \sigma\in \mathcal{S}_n
\ \ \textrm{and}\ \ \forall x_1,...,x_n\in \mathcal{N}.
\end{equation}
where $\mathcal{S}_n$ stands for the permutation group of $n$ elements.\\
The condition \eqref{NambuIdentity} may be written using skew-symmetry property  of the bracket as
\begin{eqnarray}\label{NambuIdentityLie}
  && \big[\alpha_1(x_1),....,\alpha_{n-1}(x_{n-1}),[y_1,....,y_{n}]\big]= \\ \nonumber
&& \sum_{i=1}^{n}(-1)^{i+n}\big[\alpha_1(y_1),....,\widehat{y_{i}}
  ,...,\alpha_{n-1}(y_n),[x_1,....,x_{n-1},y_i]\big],
  \end{eqnarray}
\end{definition}

In the sequel we deal sometimes with a particular class of $n$-ary Hom-Nambu algebras which we call $n$-ary multiplicative Hom-Nambu  algebras.

\begin{definition}
An \emph{$n$-ary multiplicative Hom-Nambu algebra }
(resp. \emph{$n$-ary multiplicative Hom-Nambu-Lie algebra}) is an $n$-ary Hom-Nambu algebra  (resp. $n$-ary Hom-Nambu-Lie algebra) $(\mathcal{N}, [\cdot ,..., \cdot],  \widetilde{ \alpha})$ with  $\widetilde{\alpha}=(\alpha_i)_{1\leq i\leq n-1}$
where  $\alpha_1=...=\alpha_{n-1}=\alpha$  and satisfying
\begin{equation}
\alpha([x_1,..,x_n])=[\alpha(x_1),..,\alpha(x_n)],\ \  \forall x_1,...,x_n\in \mathcal{N}.
\end{equation}
For simplicity, we will denote the $n$-ary multiplicative Hom-Nambu algebra as $(\mathcal{N}, [\cdot ,..., \cdot ],  \alpha)$ where $\alpha :\mathcal{N}\rightarrow \mathcal{N}$ is a linear map. Also by misuse of language an element  $x\in \mathcal{N}^n$ refers to  $x=(x_1,..,x_{n})$ and  $\alpha(x)$ denotes $(\alpha (x_1),...,\alpha (x_n))$.
\end{definition}

\subsection{Constructions}
In this Section we recall the construction procedures by twisting principles. The   first twisting principle, introduced for binary case  in \cite{yau}, was extend to $n$-ary case in \cite{AMS}.  The second twisting principle was introduced   in \cite{yau1}. Also we recall a construction by tensor product of symmetric totally $n$-ary Hom-associative algebra by an $n$-ary Hom-Nambu algebra given in \cite{AMS}.

The following Theorem gives a way  to construct $n$-ary multiplicative Hom-Nambu algebras  starting from a classical $n$-ary Nambu algebras  and  algebra endomorphisms.

\begin{thm}\label{twistcentral1}\cite{AMS} Let $(\mathcal{N},[\cdot,...,\cdot])$ be an $n$-ary Nambu algebra
 and let $\rho:\mathcal{N}\rightarrow \mathcal{N}$ be an n-ary Nambu
 algebra endomorphism. Then
 $(\mathcal{N},\rho\circ[\cdot ,...,\cdot ],\rho)$ is an $n$-ary multiplicative Hom-Nambu algebra.
\end{thm}
In the following we use the second twisting principal to generate new $n$-ary Hom-Nambu algebra starting from a given multiplicative $n$-ary Hom-Nambu algebra.
\begin{thm}\label{twistcentral2} Let $(\mathcal{N},[\cdot,...,\cdot],\alpha)$ be a multiplicative  $n$-ary Hom-Nambu algebra.
 Then $(\mathcal{N},\alpha^{n-1}\circ[\cdot ,...,\cdot ],\alpha^n)$, for any integer $n$,  is an $n$-ary multiplicative Hom-Nambu algebra.
\end{thm}
Now,  we define the tensor product of two n-ary Hom-algebras and prove some results involving n-ary Hom-algebras of Lie type and Hom-associative
type.

Let $A$ be a $\K$-vector space, $\mu$ be an $n$-linear map on $A$ and $\eta_i$, $i\in\{1,\cdots ,n-1\}$,  be  linear maps on $A$. A triple   $(A,\mu,\widetilde{\eta}=(\eta_1,...,\eta_{n-1}))$  is said to be a symmetric $n$-ary totally  Hom-associative algebra over $\K$ if the following identities hold
\begin{align}\label{Identityass}
&\mu(a_1,...,a_i,...,a_j,...,a_{n})=\mu(a_1,...,a_j,...,a_i,...,a_{n}),\ \forall\ i,j\in\{1,...,n\},\\
&\mu(\mu(a_1,...,a_{n}),\eta_1(a_{n+1}),...,\eta_{n-1}(a_{2n-1}))=\mu(\eta_1(a_1),\mu(a_2,...,a_{n+1}),\eta_2(a_{n+2}),...,\eta_{n-1}(a_{2n-1}))\\
&\nonumber=\cdots=\mu(\eta_1(a_1),...,\eta_{n-1}(a_{n-1}),\mu(a_n,...,a_{2n-1})),
\end{align}
where $a_1,...,a_{2n-1}\in A$.

\begin{thm}\label{prodtens}
Let $(A,\mu,\widetilde{\eta}=(\eta_1,...,\eta_{n-1}))$   be a symmetric $n$-ary totally  Hom-associative algebra and $(\mathcal{N},[\cdot,...,\cdot]_\mathcal{N},\widetilde{\alpha})$ be an $n$-ary Hom-Nambu algebras.  Then the tensor product  $A \otimes \mathcal{N}$ carries a structure  of  $n$-ary Hom-Nambu
algebra over $\K$ with respect to the  $n$-linear operation defined by
\begin{equation}
[a_1\otimes x_1 ,...,a_n\otimes x_n]=\mu(a_1,..., a_n)\otimes [x_1,...,x_n ]_\mathcal{N},\quad \text{where} \ x_l\in \mathcal{N}, a_l\in A, l\in\{1,...,n\},
\end{equation}
and linear maps $\widetilde{\zeta}=(\zeta_1,...,\zeta_{n-1})$ where $\zeta_i=\eta_i\otimes \alpha_i$, for $ i\in\{1,...,n-1\}$, defined by
\begin{equation}
\zeta_i(a\otimes x)= \eta_i(a)\otimes\alpha_i(x),\ \forall\ a\otimes x\in
A \otimes \mathcal{N}.
\end{equation}

\end{thm}

\begin{proof}

%$$\mu(\mu(a_{\sigma(1)},...,a_{\sigma(n)}),\eta_1(a_{\sigma(n+1)}),...,\eta_{n-1}(a_{\sigma(2n-1)}))
%=\mu(\mu(a_{1},...,a_{n}),\eta_1(a_{n+1}),...,\eta_{n-1}(a_{2n-1})),$$
 %for all $\sigma\in \mathcal{S}_{2n-1}$,
For $a_k\otimes x_k,\ b_l\otimes y_l\in A\otimes\mathcal{N}$, $1\leq k\leq n-1$ and $1\leq l\leq n$, we  have
\begin{align*}
&[\zeta_1(a_1\otimes x_1),...,\zeta_{n-1}(a_{n-1}\otimes x_{n-1}),[b_1\otimes y_1,...,b_{n}\otimes y_{n}]]\\
&=\mu(\eta_1(a_1),...,\eta_{n-1}(a_{n-1}),\mu( b_1,..., b_{n}))\otimes[\alpha_1( x_1),...,\alpha_{n-1}( x_{n-1}),[ y_1,...,y_{n}]_\mathcal{N}]_\mathcal{N}.
\end{align*}
The symmetry  and totally associativity of $A$ lead to
\begin{align*}
&[\zeta_1(b_1\otimes y_1),...,[a_1\otimes x_1,...,a_{n-1}\otimes x_{n-1}, b_{l}\otimes y_{l}] ,...,\zeta_{n-1}(b_{n }\otimes y_{n })]\\
&=\mu(\eta_1(b_1),...,\mu( a_1,..., a_{n-1},b_l),...,\eta_{n-1}(b_{n-1}))\otimes[[\alpha_1( y_1),...,[ x_1,...,x_{n-1},y_l]_\mathcal{N},...,\alpha_{n-1}( y_{n })]_\mathcal{N}\\
&=\mu(\eta_1(a_1),...,\eta_{n-1}(a_{n-1}),\mu( b_1,..., b_{n}))\otimes[\alpha_1( y_1),...,[ x_1,...,x_{n-1},y_l]_\mathcal{N},...,\alpha_{n-1}( y_{n })]_\mathcal{N}.
\end{align*}
Thus
\begin{align*}
&\sum_{l=1}^n[\zeta_1(b_1\otimes y_1),...,[a_1\otimes x_1,...,a_{n-1}\otimes x_{n-1}, b_{l}\otimes y_{l}] ,...,\zeta_{n-1}(b_{n }\otimes y_{n })]\\
&=\mu(\eta_1(a_1),...,\eta_{n-1}(a_{n-1}),\mu( b_1,..., b_{n}))\otimes\Big(\sum_{l=1}^n[\alpha_1( y_1),...,[ x_1,...,x_{n-1},y_l]_\mathcal{N},...,\alpha_{n-1}( y_{n })]_\mathcal{N}\Big)\\
&=\mu(\eta_1(a_1),...,\eta_{n-1}(a_{n-1}),\mu( b_1,..., b_{n}))\otimes[\alpha_1( x_1),...,\alpha_{n-1}( x_{n-1}),[ y_1,...,y_{n}]_\mathcal{N}]_\mathcal{N}.
\end{align*}
\end{proof}
\begin{cor}Let  $(A,\mu,\eta)$  be a multiplicative symmetric $n$-ary Hom-associative algebra $(\text{ i.e.}\ \eta\circ\mu=\mu\circ\eta^{\otimes n} )$ and  $(\mathcal{N},[\cdot,...,\cdot]_\mathcal{N},\widetilde{\alpha})$ be a multiplicative $n$-ary Hom-Nambu algebra.  Then $A \otimes \mathcal{N}$ is a multiplicative $n$-ary Hom-Nambu
algebra.
\end{cor}
\begin{remark}Let  $(A,\cdot)$  be a binary commutative   associative algebra  and  $(\mathcal{N},[\cdot,...,\cdot]_\mathcal{N},\widetilde{\alpha})$ be an  $n$-ary Hom-Nambu algebra.  Then the tensor product  $A \otimes \mathcal{N}$ carries a structure  of  $n$-ary Hom-Nambu
algebra over $\K$ with respect to the  $n$-linear operation defined by
\begin{equation}
[a_1\otimes x_1 ,...,a_n\otimes x_n]=(a_1\cdot...\cdot a_n)\otimes [x_1,...,x_n ]_\mathcal{N},
\end{equation}
and linear maps $\widetilde{\zeta}=(\zeta_1,...,\zeta_{n-1})$ where $\zeta_i=id\otimes \alpha_i$, for $ i\in\{1,...,n-1\}$, defined by
\begin{equation}
\zeta_i(a\otimes x)= a\otimes\alpha_i(x),\ \forall\ a\otimes x\in
A \otimes \mathcal{N}.
\end{equation}
\end{remark}
\section{Definitions and examples of  quadratic $n$-ary  Hom-Nambu algebras }
An important class of Nambu-Lie algebras, due to their appearance in a number of physical contexts, are those which possess a scalar product (here, a nondegenerate symmetric bilinear form) which is invariant under derivation.

In this section we introduce a class of Hom-Nambu-Lie algebras which possess an inner product.

\begin{definition}Let  $(\mathcal{N}, [. ,..., . ],  \widetilde{\alpha} ),\ \widetilde{\alpha}=(\alpha_i)_{1\leq i\leq n-1}$ be an $n$-ary Hom-Nambu algebra  and $B:\mathcal{N}\times\mathcal{N\rightarrow\K}$ a nondegenerate symmetric bilinear form such that, for all  $y$, $z\in \mathcal{N}$ and $x\in \mathcal{L}(\mathcal{N})$
\begin{align}\label{idinvariance}
& B([x_1,\cdots ,x_{n-1}, y] ,z)+B(y,[x_1,\cdots ,x_{n-1}, z] )=0,\\
\ & B(\alpha_i(y),z)=B(y,\alpha_i(z)),\ \ \forall i\in \{1,...,n-1\}
\end{align}
The quadruple $(\mathcal{N}, [ \cdot,...,\cdot. ], \widetilde{\alpha},B)$ is called quadratic $n$-ary Hom-Nambu algebra.
\end{definition}

\begin{remark}If $\alpha_i=Id$ for all $i\in\{1,...,n-1\}$, we recover the invariant $n$-ary Nambu algebras.
\end{remark}

\begin{definition}
An  $n$-ary Hom-Nambu algebra $(\mathcal{N}, [. ,..., . ],  \widetilde{\alpha }),\ \widetilde{\alpha}=(\alpha_i)_{1\leq i\leq n-1}$ is called \emph{Hom-quadratic} if there exists a pair $(B,\beta)$ where $B:\mathcal{N}\times\mathcal{N\rightarrow\K}$ is a  nondegenerate symmetric bilinear form  and a linear map $\beta\in End(\mathcal{N})$ satisfying
\begin{align}\label{idinvariance2}
& B([x_1,\cdots ,x_{n-1}, y] ,\beta(z))+B(\beta(y),[x_1,\cdots ,x_{n-1}, z] )=0,
\end{align}
We call \eqref{idinvariance2} the $\beta$-invariance of $B$. We recover the quadratic $n$-ary Hom-Nambu algebras when $\beta = id$. The tuple $(\mathcal{N}, [ \cdot,...,\cdot. ], \widetilde{\alpha},B,\beta)$ denotes the  Hom-quadratic $n$-ary Hom-Nambu algebra.
\end{definition}

%We provide examples of $n$-ary Hom-Nambu algebras.

\begin{example}
We consider an example ternary Hom-Nambu algebra given in \cite{yau1}.  Let $V$ be a $\mathbb{K}$-module and $B:V^{\otimes2}\rightarrow\mathbb{K}$ be a nondegenerate  symmetric bilinear form.
Suppose $\alpha: V \rightarrow V$ is an involution, that is $\alpha^2 = id$. Assume  that $B$ is $\alpha$-symmetric, that is 
$B(\alpha(x),y) =B(x,\alpha(y)) \ \  \textrm{for all}\ x, y \in V.$
We have also
$B(\alpha(x),\alpha(y))=B(\alpha^2(x),y)=B(x,y).$
Then for
any scalar $\lambda\in\mathbb{K}$, the triple product
\begin{equation}[x,y,z]_\alpha = \lambda(B(y,z)\alpha(x)-B(z,x)\alpha(y)) \ \ \ \ \ \textrm{for all} \ x,\ y,\ z \in V.\end{equation}
gives a Hom-quadratic ternary Hom-Nambu algebra $(V, [.,. , .]_\alpha,(\alpha,\alpha))$, $\alpha$-invariant by the pair $(B,\alpha)$. Indeed, for $x, y, z, t\in V$
\begin{eqnarray*}
% \nonumber to remove numbering (before each equation)
  B([x,y,z]_\alpha,\alpha(t)) &=& \lambda(B(B(y,z)\alpha(x)-B(z,x)\alpha(y),\alpha(t))) \\
    &=& \lambda(B(y,z)B(\alpha(x),\alpha(t))-B(z,x)B(\alpha(y),\alpha(t))) \\
    &=&  \lambda(B(y,z)B(x,t)-B(z,x)B(y,t)) \\
    &=& \lambda(B(\alpha(y),\alpha(z))B(x,t)-B(\alpha(z),\alpha(x))B(y,t)) \\
    &=& \lambda(B(\alpha(z),\alpha(y))B(x,t)-B(\alpha(z),\alpha(x))B(y,t)) \\
    &=& \lambda(B(\alpha(z),\alpha(y)B(x,t)-\alpha(x)By,t)))\\
    &=& -\lambda(B(\alpha(z),\alpha(x)B(y,t)-\alpha(y)B(t,x)))\\
    &=&-B(\alpha(z),[x,y,t]_\alpha).
\end{eqnarray*}

\end{example}

\begin{example}
Let $(\mathcal{N}, [\cdot,\cdot,\cdot] ,(\alpha_1, \alpha_2))$ be a $3$-dimensional  ternary  Hom-Nambu-Lie algebras,
defined with respect to a basis $\{e_1, e_2, e_3\}$ of $\mathcal{N}$ by
\begin{align}
&[e_1,e_2,e_3]= e_1+ 2e_2+ e_3,\\
&\alpha_1(e_1)=0,\ \alpha_1(e_2)=\lambda e_1+\nu e_2,\ \alpha_1(e_3)=\frac{ \lambda}{ 2}e_1+\frac{ \nu}{ 2}e_2,\\
&\alpha_2(e_1)=0,\ \alpha_2(e_2)=0,\ \alpha_2(e_3)=be_3.
\end{align}
where $\lambda,\nu,  b$ are parameters with $\lambda\nu= 0$. These $3$-dimensional admits   symmetric bilinear forms $B$ given, with respect to the previous basis,  by the following matrix
 $$M=\left(
     \begin{array}{ccc}
       -\frac{ \nu}{ \lambda} & 1 & (\frac{ \nu}{ \lambda}-2) \\
        1&-\frac{\lambda}{\nu} & (2\frac{\lambda}{\nu}-1)\\
      (\frac{ \nu}{ \lambda}-2) &(2\frac{\lambda}{\nu}-1)& (4-4\frac{\lambda}{\nu}-\frac{\nu}{\lambda}) \\
     \end{array}
   \right)
 $$
The ternary Hom-Nambu-Lie algebras $(\mathcal{N}, [\cdot,\cdot,\cdot] ,(\alpha_1, \alpha_2))$, with respect to $B$, is not quadratic because $B$ is degenerate $(det(M)= 0)$.
\end{example}
%\begin{example}
%Let $V=\mathbb{R}^4$ endowed with its Euclidean product  $B$. The $3$-bracket is defined by
%$$B([x_1,x_2,x_3],x_4)=\textrm{det}(x_1,x_2,x_3,x_4)$$
%Then  $V$ is invariant. In fact
%\begin{eqnarray*}
%% \nonumber to remove numbering (before each equation)
%  B([x_1,x_2,x_3],x_4) &=& \textrm{det}(x_1,x_2,x_3,x_4) \\
%    &=& -\textrm{det}(x_1,x_2,x_4,x_3)=-B(x_3,[x_1,x_2,x_4])
%\end{eqnarray*}
%\end{example}

\section{From quadratic $n$-ary Hom-Nambu-Lie algebra to quadratic Hom-Leibniz algebra}
In the context of  Hom-Lie algebras one
gets the class of Hom-Leibniz algebras (see \cite{MS}). A Hom-Leibniz algebra  is a triple $(V, [\cdot, \cdot], \alpha)$
consisting of a linear space $V$, a bilinear map $[\cdot, \cdot]:
V\times V \rightarrow V$ and a homomorphism $\alpha: V \rightarrow
V$ with respect to $[\cdot, \cdot]$ satisfying
\begin{equation} \label{Leibnizalgident}
[\alpha(x),[y,z]]=[[x,y],\alpha(z)]+[\alpha (y),[x,z]]
\end{equation}

We fix in the following some notations. Let $(\mathcal{N},[\cdot ,...,\cdot ],\widetilde{\alpha})$ be an $n$-ary  Hom-Nambu  algebra, we define\\
$\bullet$ a linear map  $L:\otimes^{n-1}\mathcal{N}\longrightarrow End(\mathcal{N})$ by
\begin{equation}\label{adj}L(x)\cdot z=[x_1,...,x_{n-1},z],\end{equation}
for all $x=x_1\otimes...\otimes x_{ n-1}\in\otimes^{n-1}\mathcal{N}, \ z\in \mathcal{N}$
and extending it linearly to all $\otimes^{n-1}\mathcal{N}$.  Notice that $L(x)\cdot z=ad(x)(z)$.\\

If the $n$-ary Hom-Nambu algebra $\mathcal{N}$ is multiplicative, then we define\\

$\bullet$  a linear map $\hat{\alpha}
:\otimes^{n-1}\mathcal{N}\longrightarrow\otimes^{n-1}\mathcal{N}$ by
\begin{equation}\label{mapLeibniz}\hat{\alpha}
(x)=\alpha(x_1)\otimes...\otimes\alpha(x_{n-1})\,\end{equation}

for all $x=x_1\otimes...\otimes x_{ n-1}\in\otimes^{n-1}\mathcal{N}$,\\

$\bullet$ a bilinear map $[\ ,\ ]_{\alpha}:\otimes^{n-1}\mathcal{N}\times\otimes^{n-1}\mathcal{N}\longrightarrow\otimes^{n-1}\mathcal{N}$  defined by
\begin{equation}\label{brackLei}[x ,y]_{\alpha}=L(x)\bullet_{\alpha}y=\sum_{i=0}^{n-1}\big(\alpha(y_1),...,L(x)\cdot y_i,...,\alpha(y_{n-1})\big),\end{equation}
for all $x=x_1\otimes...\otimes x_{ n-1}\in\otimes^{n-1}\mathcal{N},\ y=y_1\otimes...\otimes y_{ n-1}\in\otimes^{n-1}\mathcal{N}$\\

If $(\mathcal{N},[\cdot ,...,\cdot ],\widetilde{\alpha})$ is a multiplicative $n$-ary  Hom-Nambu-Lie  algebra, we denote by $\mathcal{L}(\mathcal{N})$ the space $\wedge^{n-1}\mathcal{N}$ and we call it  the fundamental set.

\begin{lem}\label{3.1}
Let $(\mathcal{N},[\cdot ,...,\cdot ],\widetilde{\alpha})$ be a multiplicative  $n$-ary  Hom-Nambu  algebra then the  map $L$ satisfies
\begin{equation} L([x ,y ]_{\alpha})\cdot \alpha(z)=L(\tilde{\alpha}(x))\cdot \big(L(y)\cdot z\big)-L(\tilde{\alpha}(y))\cdot \big(L(x)\cdot z\big)\end{equation}
for all $x,y\in \mathcal{L}(\mathcal{N}),\ z\in \mathcal{N}$
\end{lem}

\begin{proposition}\label{HomLeibOfHomNambu}The triple $\big(\mathcal{L}(\mathcal{N}),\ [\ ,\ ]_{\alpha},\ \hat{\alpha}
\big)$, where $[\ ,\ ]_{\alpha}$ and $\hat{\alpha}
$ are defined respectively in  \ref{mapLeibniz} and \ref{brackLei},  is a Hom-Leibniz algebra.
\end{proposition}

\begin{remark}
The invariance identity \eqref{idinvariance} of an $n$-ary Nambu algebra with respect to a bilinear form $B$ can be written
\begin{equation}\label{idinvarianceBis}
B(L(x)\cdot y ,z)+B(y,L(x)\cdot z )=0,
\end{equation}
and $\beta$-invariance identity \eqref{idinvariance2} by
\begin{equation}\label{idinvariance2bis}
B(L(x)\cdot y ,\beta(z))+B(\beta(y),L(x)\cdot z )=0,
\end{equation}
\end{remark}
\begin{proposition}
Let $(\mathcal{N},[\cdot ,...,\cdot ],\alpha,B,\alpha )$ be a Hom-quadratic multiplicative Hom-Nambu-Lie  algebra $\alpha$-invariant and  $\big(\mathcal{L}(\mathcal{N}),\ [\ ,\ ]_{\alpha},\ \hat{\alpha}
\big)$ be its associated Hom-Leibniz algebra, then the natural scalar product on $\mathcal{L}(\mathcal{N})$, $\widehat{B}$ defined by \begin{equation}
\widehat{B}(x,y)=B(x_1\wedge ...\wedge x_{n-1},y_1\wedge ...\wedge y_{n-1})=\prod_{i=0}^{n-1}B(x_i,y_i)
\end{equation}
and later extending linearly to all of $\mathcal{L}(\mathcal{N})$, is $\tilde{\alpha}$-invariant. That is, for all $x,y, z\in \mathcal{L}(\mathcal{N})$:
\begin{equation}\label{ivariance}\widehat{B}([z,x]_\alpha,\hat{\alpha}
(y))+\widehat{B}(\hat{\alpha}
(x),[z,x]_\alpha)=0.
\end{equation}
Hence, $\big(\mathcal{L}(\mathcal{N}),\ [\ ,\ ]_{\alpha},\ \hat{\alpha},\widehat{B},\hat{\alpha}
\big)$ is a Hom-quadratic Hom-Leibniz algebra.
\end{proposition}
\begin{proof}
Let $x=(x_1, ..., x_{n-1}),\ y=(y_1, ...,y_{n-1})$ and let $z\in \mathcal{L}(\mathcal{N})$. Then using equation \eqref{ivariance} we have
\begin{align*}
  \widehat{B}([z ,x ]_{\alpha},\hat{\alpha}
(y)) &= \widehat{B}(L(z)\bullet_{\alpha}x,\hat{\alpha}
(y)) \\
    &=\sum_{i=0}^{n-1}\widehat{B}(\big(\alpha(x_1),..., L(z)\cdot x_i,...,\alpha(x_{n-1})\big),\big(\alpha(y_1), ..., \alpha(y_{n-1})\big) ) \\
    &=\sum_{i=0}^{n-1}B(L(z)\cdot x_i,\alpha(y_i))\prod_{
  j=0 \
  j\neq i}^{n-1} B(\alpha(x_j),\alpha(x_j))  \\
    &= -\sum_{i=0}^{n-1}B(\alpha(x_i),L(z)\cdot y_i)\prod_{
  j=0 \
  j\neq i
}^{n-1} B(\alpha(x_j),\alpha(y_j)) \\
    &= -\sum_{i=0}^{n-1}\widehat{B}(\big(\alpha(x_1),...,\alpha(x_{n-1})\big),\big(\alpha(y_1), ...,L(z)\cdot y_i,..., \alpha(y_{n-1})\big) ) \\
    &= -\widehat{B}(\hat{\alpha}
(x),L(z)\bullet_{\alpha}y)\\& =-\widehat{B}(\hat{\alpha}
(x),[z ,y ]_{\alpha}).
\end{align*}
\end{proof}

\section{Representation theory and quadratic  $n$-ary  Hom-Nambu algebras}
In this Section we study in the general case  the  representation theory of $n$-ary Hom-Nambu  algebras introduced in \cite{AmmarSamiMakhloufNov2010} for multiplicative $n$-ary Hom-Nambu  algebras. We discuss in particular  the cases of adjoint and coadjoint representations and the connection with quadratic  $n$-ary  Hom-Nambu algebras. The results obtained in this Section generalize those given for binary case in \cite{BenayadiMakhlouf}. The representation theory of Hom-Lie algebras were independently studied in \cite{Sheng}.

\begin{definition}
A representation of an $n$-ary Hom-Nambu  algebra $(\mathcal{N},[\cdot ,...,\cdot ],\tilde{\alpha})$
on a vector space $V$ is a skew-symmetric multilinear map $\rho:\mathcal{N}^{ n-1}\longrightarrow End(V)$,
satisfying  for $x,y\in \mathcal{N}^{n-1}$ the identity
\begin{equation}\label{RepIdentity1}
\rho(\tilde{\alpha}(x))\circ\rho(y)-\rho(\tilde{\alpha}(y))\circ\rho(x)=\sum_{i=1}^{n-1}\rho(\alpha_1 (x_1),...,L(y)\cdot x_i,...,\alpha_{n-2} (x_{n-1}))\circ \nu
\end{equation}
where $\nu$ is an endomorphism on $V$.
We denote this representation by a triple $(V,\rho, \mu)$.
\end{definition}

Two representations $(V,\rho,\mu )$ and $(V',\rho',\mu' )$ of $\mathcal{N}$ are \emph{equivalent} if there exists $f:V \rightarrow V' $ an isomorphism of vector space such that $f(x\cdot v)=x\cdot ' f(v)$ and $f\circ \nu =\nu' \circ f$ where $x\cdot v=\rho(x)(v)$ and $x\cdot' v'=\rho'(x)(v')$ for $x\in \mathcal{N}^{n-1}$, $v\in V$ and $v'\in V'$. Then $V $ and $V'$ are viewed as  $\mathcal{N}^{n-1}$-modules.
\begin{example}Let $(\mathcal{N}, [\cdot  ,..., \cdot ],  \tilde{\alpha} )$ be an $n$-ary Hom-Nambu-Lie algebra. The map $ L$  defined in \eqref{adj} is a representation on $\mathcal{N}$ and  where the endomorphism  $\nu$ is the twist map $\alpha_{n-1}$. The identity \eqref{RepIdentity1} is equivalent to Hom-Nambu identity \eqref{NambuIdentity}. It is called the adjoint representation.
\end{example}

\begin{proposition}
Let $(\mathcal{N}, [.,...,], \widetilde{\alpha})$ be a $n$-ary Hom-Nambu algebra and $(V, \rho, \nu)$ be a representation of $\mathcal{N}$.
The triple $(V^*, \rho^*, \widetilde{\nu})$, where $\rho^*:\mathcal{N}^{n-1}\rightarrow End(V^*)$ is given by $\rho^*=-^t\rho$ and $\mu^*:V^*\rightarrow V^*,\ \ f\mapsto \nu^* (f)=f\circ \nu $, defines a representation of the $n$-ary Hom-Nambu-Lie algebra
$(\mathcal{N}, [.,...,], \widetilde{\alpha})$  if and only if
\begin{equation}\label{RepIdentity2}
\rho(x)\circ\rho(\tilde{\alpha}(y))-\rho(y)\circ\rho(\tilde{\alpha}(x))=\sum_{i=1}^{n-1}\nu \circ\rho(\alpha_1 (x_1),...,L(y)\cdot x_i,...,\alpha_{n-2} (x_{n-1}))
\end{equation}
\end{proposition}

\begin{proof}
Let $f\in \mathcal{N}^*$, $x, y\in \mathcal{N}^{n-1}$ and $u \in \mathcal{N}$ . We compute the right hand side of the identity \eqref{RepIdentity1}
\begin{eqnarray*}
% \nonumber to remove numbering (before each equation)
& & \rho^*(\tilde{\alpha}(x))\circ\rho^*(y)(f)(u)-\rho^*(\tilde{\alpha}(y))\circ\rho^*(x)(f)(u) \\
  & &=(\rho^*(\tilde{\alpha}(x))(\rho^*(y)(f))-\rho^*(\tilde{\alpha}(y))(\rho^*(x)(f)))(u) \\
    & &= -(\rho^*(y)(f)(\rho(\tilde{\alpha}(x)))(u))+(\rho^*(x)(f)(\rho(\tilde{\alpha}(y)))(u)) \\
    & &= f(\rho(y)(\rho(\tilde{\alpha}(x))(u)))-f(\rho(x)(\rho(\tilde{\alpha}(y))(u))) \\
    & &= f(\rho(y)(\rho(\tilde{\alpha}(x))(u))-\rho(x)(\rho(\tilde{\alpha}(y))(u))).
\end{eqnarray*}
In the other hand, the left hand side of \eqref{RepIdentity1} writes
\begin{eqnarray*}
% \nonumber to remove numbering (before each equation)
& & (\sum_{i=1}^{n-1}\rho^*(\alpha_1 (x_1),...,L(y)\cdot x_i ,...,\alpha_{n-2} (x_{n-1}))\circ\nu^*(f))(u) \\
& &= -\sum_{i=1}^{n-1}(\nu^*(f)(\rho(\alpha_1 (x_1),...,L(y)\cdot x_i ,...,\alpha_{n-2} (x_{n-1}))(u))) \\
    & &= -\sum_{i=1}^{n-1}f(\nu(\rho(\alpha_1 (x_1),...,L(y)\cdot x_i ,...,\alpha_{n-2} (x_{n-1}))(u))) \\
    & &= f(-\sum_{i=1}^{n-1}\nu(\rho(\alpha_1 (x_1),...,L(y)\cdot x_i ,...,\alpha_{n-2} (x_{n-1}))(u))).
\end{eqnarray*}
Therefore we obtain the identity \eqref{RepIdentity2}.
\end{proof}

\begin{cor}
Let $(\mathcal{N},L,\alpha_{n-1})$ be a representation of an  $n$-ary Hom-Nambu algebra $(\mathcal{N},[\cdot,...,\cdot],\widetilde{\alpha})$.
We define the map  $\widetilde{L}:\mathcal{N}^{n-1}\rightarrow End(\mathcal{N}^*)$,    for $x\in \mathcal{N}^{n-1},\ f\in \mathcal{N}^*$ and  $ y \in \mathcal{N}$,  by
($\widetilde{L}(x)\cdot f) (y)=-f(L(x)\cdot y)$.\\ Then $(\mathcal{N}^*,\widetilde{L},\alpha^* _{n-1})$ is a representation of $\mathcal{N}$ if and only if
\begin{equation}\label{RepIdentity3}
L(x)\circ L(\tilde{\alpha}( y) )-L(y)\circ L(\tilde{\alpha}( x)  )=\sum_{i=1}^{n-1}\alpha_{n-1}\circ L(\alpha_1 (x_1),...,L(y)\cdot x_i ,...,\alpha_{n-2} (x_{n-1})).
\end{equation}
\end{cor}

We establish now  a
connection between quadratic  $n$-ary Hom-Nambu algebras and representation theory.

\begin{proposition}
Let $(\mathcal{N}, [\cdot,..., \cdot],\widetilde{\alpha} )$ be an $n$-ary Hom-Nambu algebra. If there exists $B: \mathcal{N} \times \mathcal{N} \rightarrow \K $ a  bilinear form
such that the quadruple  $(\mathcal{N}, [\cdot,..., \cdot],\widetilde{\alpha} ,B)$ is a quadratic $n$-ary Hom-Nambu algebra then
\begin{enumerate}
  \item $(\mathcal{N}^*,\widetilde{L} ,\alpha^*_{n-1} )$  is a representation of $\mathcal{N}$,
  \item the representations $(\mathcal{N},L ,\alpha_{n-1} )$ and $(\mathcal{N}^*,\widetilde{L} ,\alpha^*_{n-1} )$ are isomorphic.
\end{enumerate}

\end{proposition}

\begin{proof}
To prove the first assertion, we should show that, for any $z\in \mathcal{N}$, we have
\begin{equation}\label{111}L(x)\circ L(\tilde{\alpha}(y))\cdot z -L(y)\circ L(\tilde{\alpha}(x))\cdot z =\sum_{i=1}^{n-1}\alpha_{n-1}\circ L(\alpha_1 (x_1),...,L(y)\cdot x_i ,...,\alpha_{n-2} (x_{n-1}))(z)\end{equation}
Let $u\in\mathcal{N}$
\begin{eqnarray*}
% \nonumber to remove numbering (before each equation)
    & &B(L(x)\circ L(\tilde{\alpha}(y))\cdot z -L(y)\circ L(\tilde{\alpha}(x))\cdot z ,u)  \\
    & &=B(L(x)\circ L(\tilde{\alpha}(y))\cdot z ,u)-(L(y)\circ L(\tilde{\alpha}(x))\cdot z ,u)  \\
      & &=B( L(\tilde{\alpha}(y))\cdot z ,L(x)\cdot u )-( L(\tilde{\alpha}(x))\cdot z ,L(y)\cdot u )  \\
      & &=B( z,L(\tilde{\alpha}(y))\circ L(x)\cdot u )-( z,L(\tilde{\alpha}(x))\circ L(y)\cdot u )  \\
    & &=B( z,L(\tilde{\alpha}(y))\circ L(x)(u)-L(\tilde{\alpha}(x))\circ L(y)(u))
\end{eqnarray*}
and
\begin{eqnarray*}
% \nonumber to remove numbering (before each equation)
    & &B(\alpha_{n-1}\circ L(\alpha_1 (x_1),...,L(y)\cdot x_i ,...,\alpha_{n-2} (x_{n-1}))(z),u)  \\
    & &B(\alpha_{n-1}\circ L(\alpha_1 (x_1),...,L(y)\cdot x_i ,...,\alpha_{n-2} (x_{n-1}))(z),u)  \\
      & &=B(L(\alpha_1 (x_1),...,L(y)\cdot x_i ,...,\alpha_{n-2} (x_{n-1}))(z),\alpha_{n-1}(u))  \\
      & &=B(z,L(\alpha_1 (x_1),...,L(y)\cdot x_i ,...,\alpha_{n-2} (x_{n-1}))\circ\alpha_{n-1}(u))  \\
    & &=B( z,L(\alpha_1 (x_1),...,L(y)\cdot x_i ,...,\alpha_{n-2} (x_{n-1}))\circ\alpha_{n-1}(u))
\end{eqnarray*}
Since $B$ is bilinear, then
\begin{eqnarray*}
% \nonumber to remove numbering (before each equation)
    & & B(\sum_{i=1}^{n-1}\alpha_{n-1}\circ L(\alpha_1 (x_1),...,L(y)\cdot x_i ,...,\alpha_{n-2} (x_{n-1}))(z),u) \\
    & & =B( z,\sum_{i=1}^{n-1}L(\alpha_1 (x_1),...,L(y)\cdot x_i ,...,\alpha_{n-2} (x_{n-1}))\circ\alpha_{n-1}(u))
\end{eqnarray*}

Hence
\begin{eqnarray*}
% \nonumber to remove numbering (before each equation)
    & & B(L(x)\circ L(\tilde{\alpha}(y))\cdot z -L(y)\circ L(\tilde{\alpha}(x))\cdot z -\sum_{i=1}^{n-1}\alpha_{n-1}\circ L(\alpha_1 (x_1),...,L(y)\cdot x_i ,...,\alpha_{n-2} (x_{n-1}))(z),u) \\
    & & =B(z,L(\tilde{\alpha}(y))\circ L(x)\cdot u -L(\tilde{\alpha}(x))\circ L(y)\cdot u -\sum_{i=1}^{n-1}L(\alpha_1 (x_1),...,L(y)\cdot x_i ,...,\alpha_{n-2} (x_{n-1}))\circ\alpha_{n-1}(u))\\
    & &=0
\end{eqnarray*}
Since $B$ is nondegenerate then we have the identity \eqref{111} .\\
For the second assertion we consider the map $\psi : \mathcal{N} \rightarrow \mathcal{N}^*$ defined by $x \mapsto B(x, \cdot)$ which is bijective since
B is nondegenerate and prove that it is also a module morphism.
\end{proof}

\section{ Constructions of Quadratic $n$-ary  Hom-Nambu algebras }

\subsection{Twisting principles}
We show in the following some constructions  of Hom-quadratic Hom-Nambu-Lie algebras starting from an ordinary Nambu-Lie algebra and from
tensor product of Hom-quadratic commutative Hom-associative algebra and Hom-quadratic  Hom-Nambu-Lie algebra considered in Theorem \ref{prodtens}.\\

Let $(\mathcal{N},[\cdot,...,\cdot],B)$ be a quadratic  $n$-ary Nambu  algebras. We denote $Aut_S(\mathcal{N},B)$  by the set of symmetric automorphisms
of $\mathcal{N}$ with respect of $B$, that is automorphisms $f :\mathcal{N}\rightarrow\mathcal{N} $ such that $B(f(x), y) = B(x, f(y)), \forall x, y \in\mathcal{N} $.
\begin{proposition}\label{quadmultip}
Let $( \mathcal{N},[\cdot ,...,\cdot ]
,B)$ be a quadratic Nambu  algebra and $\rho\in Aut_S(\mathcal{N},B)$. 

Then   $(\mathcal{N},[\cdot ,...,\cdot ]
_\rho,\tilde{\rho},B,\rho)$ where
\begin{align}
&[\cdot ,...,\cdot ]
_\rho=\rho\circ[\cdot ,...,\cdot ]
\end{align}
 is  Hom-quadratic Hom-Nambu algebra,  and $(\mathcal{N},[\cdot ,...,\cdot ]
_\rho,\tilde{\rho},B_\rho)$ where
\begin{align}
&B_\rho(x,y)=B(\rho(x),y)
\end{align}is a quadratic Hom-Nambu algebra.
\end{proposition}

\begin{proof}
Let $x=(x_1, ...,x_{n-1})\in \mathcal{N}^{\otimes n-1}$ et $y_1,y_2\in \mathcal{N}$,
\begin{align*}
  B([x_1,...,x_{n-1},y_1]_{\rho},\rho(y_2)) &= B([\rho(x_1),...,\rho(x_{n-1)},\rho(y_1)],\rho(y_2))  \\
    &= -B(\rho(y_1),[\rho(x_1),...,\rho(x_{n-1}),\rho(y_2)])  \\
    &= -B(\rho(y_1),\rho\circ[x_1,...,x_{n-1},y_2])  \\
    &= -B(\rho(y_1),[x_1,...,x_{n-1},y_2]_{\rho}).
\end{align*}
In the other hand we have
\begin{align*}
  B_\rho([x_1,...,x_{n-1},y_1]_{\rho},y_2) &= B(\rho[\rho(x_1),...,\rho(x_{n-1)},\rho(y_1)],y_2)  \\
  &=B([\rho(x_1),...,\rho(x_{n-1)},\rho(y_1)],\rho(y_2))\\
    &= -B(\rho(y_1),[\rho(x_1),...,\rho(x_{n-1}),\rho(y_2)])  \\
    &= -B(\rho(y_1),\rho\circ[x_1,...,x_{n-1},y_2])  \\
    &= -B_\rho(y_1,[x_1,...,x_{n-1},y_2]_{\rho}).
\end{align*}
\end{proof}

\begin{proposition}
Let $( \mathcal{N},[\cdot ,...,\cdot ]
,\alpha,B)$ be a quadratic multiplicative Hom-Nambu  algebra. Then the  $(\mathcal{N},\alpha^{n-1}\circ[\cdot ,...,\cdot ]
,\alpha^n,B,\alpha^{n-1})$ is  Hom-quadratic Hom-Nambu algebra  and $(\mathcal{N},\alpha^{n-1}\circ[\cdot ,...,\cdot ]
,\alpha^n,B_\alpha)$, where
\begin{equation}
B_\alpha(x,y)=B(\alpha^{n-1}(x),y)=B(x,\alpha^{n-1}(y)),
\end{equation}
 is quadratic Hom-Nambu algebra,

\end{proposition}
\begin{proof}
Using the second twisting centrale construction in Theorem \ref{twistcentral2}, $(\mathcal{N},\alpha^{n-1}\circ[\cdot ,...,\cdot ]
,\alpha^n)$ is a  Hom-Nambu algebra. Let now $x_i,y,z\in\mathcal{N}$, $i\in\{1,...,n-1\}$, we have
$$B_\alpha(\alpha^{n-1}(y),z)=B(\alpha^{2n-2}(y),z)=B(\alpha^{n-1}(y),\alpha^{n-1}(z))=B_\alpha(y,\alpha^{n-1}(z)).$$
In the other hands, we have
\begin{align*}
B_\alpha(\alpha^{n-1}\circ[x_1,...,x_{n-1},y],z)&=B(\alpha^{n-1}\circ[x_1,...,x_{n-1},y],\alpha^{n-1}(z))\\
&=B([\alpha^{n-1}(x_1),...,\alpha^{n-1}(x_{n-1}),\alpha^{n-1}(y)],\alpha^{n-1}(z))\\
&=-B(\alpha^{n-1}(y),[\alpha^{n-1}(x_1),...,\alpha^{n-1}(x_{n-1}),\alpha^{n-1}(z)])\\
&=-B(\alpha^{n-1}(y),\alpha^{n-1}\circ[x_1,...,x_{n-1},z)]\\
&=-B_\alpha(y,\alpha^{n-1}\circ[x_1,...,x_{n-1},z)].
\end{align*}
Therefore $B_{\alpha}$ is invariant.
\end{proof}
\subsection{$T^*$-Extension of $n$-ary Hom-Nambu algebras}
We provide here a construction of Hom-Nambu algebra $\mathcal{L}$ which is a generalization of the trivial $T^*$-extension
introduced in \cite{Bord,medina}.

\begin{thm}
Let $(\mathcal{N}, [\cdot,... ,\cdot ]_\mathcal{N},B)$ be a quadratic  $n$-ary Nambu-Lie algebra and $\mathcal{N}^*$ be the underlying dual vector space. The vector
space $\mathcal{L} = \mathcal{N}\oplus\mathcal{N}^*$ equipped with the following product $[\cdot,...,\cdot]_\mathcal{L}:\mathcal{L}^n\rightarrow\mathcal{L}$ given, for $u_i=x_i+f_i\in\mathcal{L}$ where $i\in\{1,...,n\}$ by
\begin{equation}
[u_1,...,u_n]_\mathcal{L}=[x_1,...,x_n]_\mathcal{N}+\sum_{i=1}^n(-1)^{i+n+1}f_i\circ L(x_1,...,\widehat{x_i},...,x_n),
\end{equation}
and a bilinear form
\begin{equation}
B_\mathcal{L}:\begin{array}{c}
    \mathcal{L}\times\mathcal{L}\longrightarrow\mathcal{L} \\
    B_\mathcal{L}(x+f,y+g)=B(x,y)+f(y)+g(x)
  \end{array}
\end{equation}
is a quadratic $n$-ary Nambu algebra.
\end{thm}

\begin{proof}
$$$$
$\star)$ Set $u_i=x_i+f_i\in\mathcal{L}$ and $v_i=y_k+g_k\in\mathcal{L}$. We show the following  Nambu identity on $\mathcal{L}$ 
\begin{align}\label{fondIdentity}
&[u_1,...,u_{n-1},[v_1,...,v_n]_\mathcal{L}]_\mathcal{L}
=\sum_{l=1}^{n}(-1)^{l+n}[v_1,...,\widehat{v_l},...,v_{n},[u_1,...,u_{n-1},v_l]_\mathcal{L}]_\mathcal{L}.
\end{align}
Let us compute first $[u_1,...,u_{n-1},[v_1,...,v_n]_\mathcal{L}]_\mathcal{L}$. This is given by
\begin{align*}
&[u_1,...,u_{n-1},[v_1,...,v_n]_\mathcal{L}]_\mathcal{L}\\&=[x_1,...,x_{n-1},[y_1,...,y_n]_\mathcal{N}]_\mathcal{N}+\sum_{i=1}^{n-1}(-1)^{i+n+1}f_i\circ L(x_1,...,\widehat{x_i},...,x_{n-1},[y_1,...,y_n]_\mathcal{N})\\
&+\sum_{i=1}^n(-1)^{i+n}g_i\circ L(y_1,...,\widehat{y_i},...,y_n)\circ L(x_1,...,x_{n-1}).
\end{align*}
Hence the  right hand side \eqref{fondIdentity} give, for any $l\in\{1,...,n\}$
\begin{align*}
&[v_1,...,\widehat{v_l},...,v_{n},[u_1,...,u_{n-1},v_l]_\mathcal{L}]_\mathcal{L}
=[y_1,...,\widehat{y_l},...,y_{n},[x_1,...,x_{n-1},y_l]_\mathcal{N}]_\mathcal{N}\\&+\sum_{ i=1 }^{n-1}(-1)^{i+n}f_i\circ L(x_1,...,\widehat{x_i},...,x_{n-1},y_l)\circ L(y_1,...,\widehat{y_l},...,y_{n})\\
&+\sum_{i=1}^n(-1)^{i+n+1}g_i\circ L(y_1,...,\widehat{y_i},...,\widehat{y_l},...,y_{n},[x_1,...,x_{n-1},y_l]_\mathcal{N})\\
        &+g_l\circ L(x_1,...,x_{n-1})\circ L(y_1,...,\widehat{y_l},...,y_{n}).
\end{align*}

$$$$

$\star)$ Using the Nambu identity on $\mathcal{N}$, we have
\begin{align*}
&[x_1,...,x_{n-1},[y_1,...,\widehat{y_l},...,y_n,z]_\mathcal{N}]_\mathcal{N}=\sum_{i=1,i\neq l}^n(-1)^{i+n}[y_1,...,\widehat{y_i},...,\widehat{y_l},...,y_{n-1},[x_1,...,x_{n-1},y_i]_\mathcal{N},z]_\mathcal{N}\\
&+[y_1,...,\widehat{y_l},...,y_{n-1},[x_1,...,x_{n-1},z]_\mathcal{N}]_\mathcal{N}.
\end{align*}
Equivalently 
\begin{align*}
&L(y_1,...,\widehat{y_l},...,y_{n})\circ L(x_1,...,x_{n-1})\\
&=L(x_1,...,x_{n-1})\circ L(y_1,...,\widehat{y_l},...,y_{n})+\sum_{i=1}^{n-1}(-1)^{i+n+1}L(y_1,...,\widehat{y_i},...,\widehat{y_l},...,y_{n},[x_1,...,x_{n-1},y_l]).
\end{align*}
Thus for any $l\in\{1,...n\}$
\begin{align*}
&g_l\circ L(y_1,...,\widehat{y_l},...,y_{n})\circ L(x_1,...,x_{n-1})-g_l\circ L(x_1,...,x_{n-1})\circ L(y_1,...,\widehat{y_l},...,y_{n})\\
&=\sum_{i=1}^{n-1}(-1)^{i+n}g_l\circ L(y_1,...,\widehat{y_i},...,\widehat{y_l},...,y_{n},[x_1,...,x_{n-1},y_l]).
\end{align*}
$$$$
$\star)$ In the other hand we show that, for  $k\in\{1,...n\}$
\begin{align*}
&-f_k\circ L(x_1,...,\widehat{x_k},...,x_{n-1},[y_1,...,y_n]_\mathcal{N})\\
&=\sum_{i=1}^n(-1)^{i+n}f_k\circ L(x_1,...,\widehat{x_k},...,x_{n-1},y_i)\circ L(y_1,...,\widehat{y_i},...,y_{n}).
\end{align*}

Using the Namby identity \eqref{fondIdentity} on $\mathcal{N}$ and the invariance of $B$, we obtain
\begin{align*}
&B([x_1,...,\widehat{x_k},...,x_{n},[y_1,...,y_n]_\mathcal{N}]_\mathcal{N},z)=
B(x_{n},[x_1,...,\widehat{x_k},...,x_{n-1},[y_1,...,y_n]_\mathcal{N},z]_\mathcal{N}).
\end{align*}
Hence
\begin{align*}
&B(\sum_{i=1}^n(-1)^{i+n}[y_1,...,\widehat{y_i},...,y_{n},[x_1,...,\widehat{x_k},...,x_{n},y_i]_\mathcal{N}]_\mathcal{N},z)\\
&=\sum_{i=1}^n(-1)^{i+n}B([y_1,...,\widehat{y_i},...,y_{n},[x_1,...,\widehat{x_k},...,x_{n},y_i]_\mathcal{N}]_\mathcal{N},z)\\
&=-\sum_{i=1}^n(-1)^{i+n}B([x_1,...,\widehat{x_k},...,x_{n},y_i]_\mathcal{N},[y_1,...,\widehat{y_i},...,y_{n},z]_\mathcal{N})\\
&=-\sum_{i=1}^n(-1)^{i+n}B(x_{n},[x_1,...,\widehat{x_k},...,x_{n-1},y_i,[y_1,...,\widehat{y_i},...,y_{n},z]_\mathcal{N}]_\mathcal{N})\\
&=-B(x_{n},\sum_{i=1}^n(-1)^{i+n}[x_1,...,\widehat{x_k},...,x_{n-1},y_i,[y_1,...,\widehat{y_i},...,y_{n},z]_\mathcal{N}]_\mathcal{N}).
\end{align*}
Since $B$ is nondegenerate, then
\begin{align*}
&-[x_1,...,\widehat{x_k},...,x_{n-1},[y_1,...,y_n]_\mathcal{N},z]_\mathcal{N}\\
&=\sum_{i=1}^n(-1)^{i+n}
[x_1,...,\widehat{x_k},...,x_{n-1},y_i,[y_1,...,\widehat{y_i},...,y_{n},z]_\mathcal{N}]_\mathcal{N},
\end{align*}

and equivalently 
$$-L(x_1,...,\widehat{x_k},...,x_{n-1},[y_1,...,y_n]_\mathcal{N})=\sum_{i=1}^n(-1)^{i+n} L(x_1,...,\widehat{x_k},...,x_{n-1},y_i)\circ L(y_1,...,\widehat{y_i},...,y_{n}).$$
Finally,  the Nambu identity \eqref{fondIdentity} is satisfied. Thus $(\mathcal{L},[\cdot,...,\cdot]_\mathcal{L})$ is an $n$-ary  Nambu algebra.
\end{proof}

\begin{thm}
Let $(\mathcal{N},[\cdot,...,\cdot]_\mathcal{N},B)$ be a quadratic $n$-ary Nambu-Lie algebra where $\alpha\in Aut_S(\mathcal{N},B)$  is an involution. Then
  $(\mathcal{L},[\cdot,...,\cdot]_\Omega,\widetilde{\Omega},B_\mathcal{L},\Omega)$, where 
  $
  \Omega:\mathcal{L}\rightarrow\mathcal{L}$, 
 $ x+f\rightarrow \Omega(x+f)=\alpha(x)+f\circ\alpha
$
   and
  $[\cdot,...,\cdot]_\Omega=\Omega\circ[\cdot,...,\cdot]_\mathcal{L}$,  is a Hom-quadratic multiplicative $n$-ary   Hom-Nambu  algebra.
  
\end{thm}

\begin{proof}Let $x_1,...,x_n\in\mathcal{N}$ and $f_1,...,f_n\in\mathcal{N}^*$,
\begin{align*}
&\Omega[x_1+f_1,...,x_n+f_n]_\mathcal{L}=\alpha[x_1,...,x_n]_\mathcal{N}+\sum_{i=1}^n(-1)^if_i\circ L(x_1,...,\widehat{x_i},...,x_n)\circ\alpha,\\
&[\Omega (x_1+f_1),...,\Omega(x_n+f_n)]_\mathcal{L}=[\alpha(x_1),...,\alpha(x_n)]_\mathcal{N}+\sum_{i=1}^n(-1)^if_i\circ\alpha\circ L(\alpha(x_1),...,\widehat{x_i},...,\alpha(x_n)).
\end{align*}
That is for all $z\in\mathcal{N}$
\begin{align*}
\alpha\circ L(\alpha(x_1),...,\widehat{x_i},...,\alpha(x_n))(z)&=\alpha[\alpha(x_1),...,\widehat{x_i},...,\alpha(x_n),z]_\mathcal{N}\\
&=[\alpha^2(x_1),...,\widehat{x_i},...,\alpha(x_n)^2,\alpha(z)]_\mathcal{N}\\
&=[x_1,...,\widehat{x_i},...,x_n,\alpha(z)]_\mathcal{N}\\
& =L(x_1,...,\widehat{x_i},...,x_n)\cdot\alpha(z).
\end{align*}
Then $\Omega[x_1+f_1,...,x_n+f_n]_\mathcal{L}=[\Omega (x_1+f_1),...,\Omega(x_n+f_n)]_\mathcal{L}$.\\
In the following we show that $\Omega$ is symmetric with respect to $B_\mathcal{L}$. Indeed, let $x, y\in\mathcal{N}$ and $f, h\in\mathcal{N}$ \begin{align*}
B_\mathcal{L}(\Omega(x+f),y+h)&=B_\mathcal{L}(\alpha(x)+f\circ\alpha,y+h)\\
&=B(\alpha(x),y)+f\circ\alpha(y)+h\circ\alpha(x)\\
&=B(x,\alpha(y))+f\circ\alpha(y)+h\circ\alpha(x)\\
&=B_\mathcal{L}(x+f,\alpha(y)+h\circ\alpha)=B_\mathcal{L}(x+f,\Omega(y+h))
\end{align*}
Thus, using Proposition \ref{quadmultip}, $(\mathcal{L},[\cdot,...,\cdot]_\Omega,\widetilde{\Omega},B_\mathcal{L},\Omega)$  is a Hom-quadratic multiplicative   $n$-ary Hom-Nambu  algebra. We have also that  $(\mathcal{L},[\cdot,...,\cdot]_\Omega,\Omega,B_{\mathcal{L},\Omega})$, where $B_{\mathcal{L},\Omega}(u,v)=B_{\mathcal{L}}(\Omega(u),v)$, for all $u,v\in\mathcal{L}$,  is a quadratic multiplicative  $n$-ary Hom-Nambu  algebra.

\end{proof}

\subsection{Tensor product}
Let $(A,\mu,\widetilde{\eta},B_A,\beta_A)$ be a Hom-quadratic symmetric $n$-ary totally Hom-associative algebra, that is a symmetric $n$-ary totally Hom-associative algebra together with a symmetric nondegenerate form satisfying  the following assertions
\begin{align}
&B_A(\eta_i(a),b)=B_A(a,\eta_i(b)),\ \textrm{for\ all} \ i\in\{1,...,n-1\}\\
&B_A(\mu(a_1,...a_{n-1}, b),\beta_A(c))=B_A(\beta_A(a),\mu(a_1,...a_{n-1}, c)),\ \textrm{for\ all}\ a_i,b,c\in A, \ i\in\{1,...,n-1\}.
\end{align}
We discuss now the tensor product as in Proposition \ref{prodtens}.
\begin{thm} Let $(\mathcal{N},[\cdot,...,\cdot]_\mathcal{N},\alpha,B_\mathcal{N},\beta_\mathcal{N})$ be a Hom-quadratic $n$-ary Hom-Nambu algebra, then $(A\otimes\mathcal{N},[\cdot,...,\cdot],\zeta,\widetilde{B},\omega)$  where 
\begin{align}
&\widetilde{B}(a\otimes x ,b\otimes y)=B_A(a,b)B_\mathcal{N}(x,y),\\
&\omega(a\otimes x)=\beta_A(a)\otimes\beta_\mathcal{N}(x).
\end{align}
is a Hom-quadratic $n$-ary Hom-Nambu algebra.
\end{thm}

\subsection{Hom-quadratic Hom-Nambu-Lie algebras induced by Hom-quadratic Hom-Lie algebras }In  \cite{jas} the authors provided a construction procedure of ternary Hom-Nambu-Lie algebras starting from a bilinear bracket of a Hom-Lie algebra and a trace function satisfying certain compatibility conditions including the twisting map.\\

The aim of this section is to prove that this procedure is still true for quadratic Hom-Nambu-Lie algebra.\\
First we recall the result in \cite{jas}.
\begin{definition}
Let $ (V,[\cdot , \cdot ])$
be a binary algebra and let $\tau:V\rightarrow\mathbb{K}$ be a linear form. The trilinear map $[\cdot , \cdot ,\cdot]_\tau:V\times V\times V\rightarrow V$ is defined as
\begin{equation}
[x,y,z]_\tau=\tau(x)[y,z]+\tau(y)[z,x]+\tau(z)[x,y].
\end{equation}
\end{definition}
\begin{remark}If the bilinear multiplication $[\cdot , \cdot]$  is skew-symmetric, then the trilinear map $[\cdot , \cdot ,\cdot]_\tau$ is skew-symmetric as well.
\end{remark}
\begin{thm}\cite{jas} Let  $(V, [\cdot , \cdot],\alpha)$ be a Hom-Lie algebra and $\gamma:V\longrightarrow V$ be a linear map. Furthermore, assume that $\tau$ is a trace function on $V$ fulfilling
\begin{eqnarray}
% \nonumber to remove numbering (before each equation)
  \tau(\alpha(x))\tau(y) &=& \tau(x)\tau(\alpha(y)), \\
  \tau(\gamma(x))\tau(y) &=& \tau(x)\tau(\gamma(y)), \\
  \tau(\alpha(x))\gamma(y) &=& \tau(\gamma(x))\alpha(y),
\end{eqnarray}
for all $x,y\in V$. Then $(V,[\cdot ,\cdot, \cdot]_\tau,(\alpha,\gamma))$ is a ternary Hom-Nambu-Lie algebra, and we say that it induced by $(V, [\cdot , \cdot],\alpha)$.
\end{thm}
\begin{proposition}
Let $(V, [\cdot , \cdot ],\alpha,B,\beta)$ be a Hom-quadratic Hom-Lie algebra  satisfying
\begin{align}
&B(\alpha(x),y)=B(x,\alpha(y)),\\
&B(\gamma(x),y)=B(x,\gamma(y)),\\
&\tau(x)B(\beta(y),z)-\tau(y)B(\beta(x),z) =0\ \ \   \textrm{for all } x,y,z\in V.
\end{align}
Then  $(V,[\cdot , \cdot,\cdot ]_\tau,(\alpha,\gamma),B,\beta)$ is a Hom-quadratic ternary  Hom-Nambu-Lie algebra.
\end{proposition}
\begin{proof}
Let $x_1,x_2, y_1,y_2 \in V$
\begin{eqnarray*}
% \nonumber to remove numbering (before each equation)
  B([x_1,x_2,y_1]_\tau,\beta(y_2)) &=& \tau(x_1)B([x_2,y_1],\beta(y_2))-\tau(x_2)B([x_1,y_1],\beta(y_2)) \\
    &+& \tau(y_1)B([x_1,x_2],\beta(y_2)).
\end{eqnarray*}
\begin{eqnarray*}
% \nonumber to remove numbering (before each equation)
  B(\beta(y_1),[x_1,x_2,y_2]_\tau) &=& \tau(x_1)B(\beta(y_1),[x_2,y_2])-\tau(x_2)B(\beta(y_1),[x_1,y_2])  \\
    &-& \tau(y_2)B(\beta(y_1),[x_1,x_2]).
\end{eqnarray*}
Since $B$ is symmetric, then
$$B([x_1,x_2,y_1]_\tau,\beta(y_2))+B(\beta(y_1),[x_1,x_2,y_2]_\tau)=0. $$
\end{proof}
\subsection{Quadratic  Hom-Nambu algebras of higher arities}
The purpose of this section is to observe that every Hom-quadratic multiplicative $n$-ary Hom-Nambu algebra
gives rise to a sequence of quadratic multiplicative Hom-Nambu algebras of increasingly higher
arities. The construction of this sequence was given first in \cite{yau1}.
\begin{thm}Let $(\mathcal{N}, [\cdot, . . ., \cdot], \alpha,B,\beta)$ be a Hom-quadratic multiplicative $n$-ary Hom-Nambu algebra. Define the $(2n-1)$-
ary product
\begin{equation}
[x_1,...,x_{2n-1}]^{(1)}=[[x_1,...,x_n],\alpha(x_{n+1}),...,\alpha(x_{2n-1})]\ \ \ \ \textrm{for}\ \ x_i\in V
\end{equation}
Then $\mathcal{N}^1=(\mathcal{N}, [\cdot, . . . ,\cdot ]^{(1)}, \alpha^2,B,\beta')$, where $\beta'=\beta\alpha$, is a Hom-quadratic multiplicative $(2n-1)$-ary Hom-Nambu algebra.
\end{thm}
\begin{proof}
For the proof of  the $(2n - 1$)-ary Hom-Nambu identity and the multiplicativity for $\mathcal{N}^1$, see \cite{yau1}.\\
Let $x_1,...,x_{2n-2},y_1,y_2\in \mathcal{N}$
\begin{eqnarray*}
% \nonumber to remove numbering (before each equation)
  B([x_1,...,x_{2n-2},y_1]^{(1)},\beta'(y_2)) &=& B([[x_1,...,x_n],\alpha(x_{n+1}),...,\alpha(x_{2n-2}),\alpha(y_1)],\beta(\alpha(y_2))) \\
    &=&- B(\beta(\alpha(y_1)),[[x_1,...,x_n],\alpha(x_{n+1}),...,\alpha(x_{2n-2}),\alpha(y_2)]) \\
    &=&- B(\beta'(y_1),[x_1,...,x_{2n-2},y_2]^{(1)}).
\end{eqnarray*}
Hence
$$B([x_1,...,x_{2n-2},y_1]^{(1)},\beta'(y_2)) +B(\beta'(y_1),[x_1,...,x_{2n-2},y_2]^{(1)})=0.$$
\end{proof}
\begin{cor}
Let $(\mathcal{N}, [\cdot, . . . ,\cdot ], \alpha,B,\beta)$ be a Hom-quadratic multiplicative $n$-ary Hom-Nambu algebra. For $k\geq 1$ define the $(2^k(n-1)+1)$-
ary product $[\cdot, . . . ,\cdot ]^{(k)}$ inductively by setting $[\cdot, . . . ,\cdot ]^{(0)} = [\cdot, . . . , \cdot]$ and
\begin{equation}
[x_1,...,x_{2^k(n-1)+1}]^{(k)}=[[x_1,...,x_{2^{k-1}(n-1)+1}]^{(k-1)},\alpha^{2^{k-1}}(x_{2^{k-1}(n-1)+2}),...,\alpha^{2^{k-1}}(x_{2^k(n-1)+1})]^{(k-1)}
\end{equation}
for all $x_i\in \mathcal{N}$.\\Then $\mathcal{N}^k=(\mathcal{N}, [\cdot, . . . , \cdot]^{(k)}, \alpha^{2^k},B,\beta')$, where $\beta'=\beta\alpha^{2^{k-1}}$, is a Hom-quadratic multiplicative $(2^k(n-1)+1)$-ary Hom-Nambu algebra. 

\end{cor}
\subsection{Quadratic  Hom-Nambu algebras of lower arities}
The purpose of this section is to observe that, under suitable assumptions, a quadratic $n$-ary Hom-Nambu algebra with $n \geq 3$ reduces to an quadratic $(n - 1)$-ary Hom-Nambu algebra. We use the construction given in  \cite{yau1}.
\begin{thm}
Let $n \geq 3$ and $(\mathcal{N}, [\cdot, . . . ,\cdot ], \alpha = (\alpha_1, . . . , \alpha_{n-1}),B,\beta)$ be a Hom-quadratic $n$-ary Hom-Nambu algebra.
Suppose $a \in \mathcal{N}$ satisfies
$$\alpha_1(a)=a\ \ \textrm{and}\ \ [a,x_1,...,x_{n-2},a]=0\ \ \ \textrm{for\ all}\ \ x_i\in \mathcal{N}$$
Then
$\mathcal{N}_a=(\mathcal{N}, [\cdot, . . . , \cdot]_a, \alpha_a = (\alpha_2, . . . , \alpha_{n-1}),B,\beta),$  where
$$[x_1,...,x_{n-1}]_a=[a,x_1,...,x_{n-1}]\ \ \ \textrm{for\ all}\ \ x_i\in \mathcal{N},$$
is a Hom-quadratic $(n - 1)$-ary Hom-Nambu algebra.
\end{thm}
\begin{proof}
Using \cite{yau1}, $\mathcal{N}_a=(\mathcal{N}, [\cdot , . . . ,\cdot ]_a, \alpha_a = (\alpha_2, . . . , \alpha_{n-1}))$
is an  $(n - 1)$-ary Hom-Nambu algebra.\\
Let $x_1,...,x_{n-2},y_1,y_2\in \mathcal{N}$, then
\begin{eqnarray*}
% \nonumber to remove numbering (before each equation)
  B([x_1,...,x_{n-2},y_1]_a,\beta(y_2)) &=& B([a,x_1,...,x_{n-2},y_1],\beta(y_2)) \\
    &=& -B(\beta(y_1),[a,x_1,...,x_{n-2},y_2]) \\
    &=& -B(\beta(y_1),[x_1,...,x_{n-2},y_2]_a).
\end{eqnarray*}
Hence
$$B([x_1,...,x_{n-2},y_1]_a,\beta(y_2)) +B(\beta(y_1),[x_1,...,x_{n-2},y_2]_a)=0.$$
\end{proof}
\begin{cor}
Let $(\mathcal{N}, [\cdot, . . . , \cdot], \alpha = (\alpha_1, . . . , \alpha_{n-1}),B,\beta)$ be a Hom-quadratic $n$-ary Hom-Nambu algebra, with $n \geq 3$.
Suppose for some $k\in\{1,...,n-2\}$ there exist  $a_i \in L$ for $1\leq i\leq k$ satisfying
$$\alpha_i(a_i)=a_i\ \ \textrm{for} \ \ 1\leq i\leq k$$
and
$$[a_1,...,a_j,x_{j+1},...,x_{n-1},a_j]=0\ \ \ \textrm{for}\ \ 1\leq j\leq k \ \ \textrm{and\ all}\ \ x_l\in \mathcal{N}$$
Then
$\mathcal{N}_k=(\mathcal{N}, [\cdot, . . . ,\cdot ]_k, \alpha_ k= (\alpha_{k+1}, . . . , \alpha_{n-1}),B,\beta)$, where
$$[x_1,...,x_{n-1}]_k=[a_1,...,a_k,x_k,...,x_{n-1}]\ \ \ \textrm{for\ all}\ \ x_i\in \mathcal{N}$$
is a Hom-quadratic $(n - k)$-ary Hom-Nambu algebra.
\end{cor}

\subsection{ Ternary Nambu algebras arising from the real Faulkner construction}
Let $(\mathfrak{g},[\cdot,\cdot ],B)$ be a real finite-dimensional quadratic  Lie algebra and let $\mathfrak{g}^*$ be the  dual of $\mathfrak{g}$.
We denote by $\langle-,-\rangle$ the dual pairing between $\mathfrak{g}$ and $\mathfrak{g}^*$.\\
For all $x\in\mathfrak{g}$ and $f\in\mathfrak{g}^*$ we define
an element $\phi(x\otimes f)\in\mathfrak{g}$ by
\begin{equation}\label{sc}B(y,\phi(x\otimes f))=\langle[y,x],f\rangle=f([y,x])\ \textrm{for\ all} \ y\in\mathfrak{g}.
\end{equation}
Extending $\phi$ linearly, defines a $\mathfrak{g}$-equivariant map
$\phi : \mathfrak{g} \otimes \mathfrak{g}^*\longrightarrow \mathfrak{g}$, which is surjective. To
lighten the notation we will write $ \phi(x, f)$ for $\phi(x \otimes f)$ in the sequel. The $\mathfrak{g}$-equivariance of $\phi$
is equivalent to
\begin{equation}\label{aaa}[\phi(x, f),\phi(y, g)]=\phi([\phi(x, f), y], g)+\phi(y, \phi(x, f)\cdot g),
\end{equation}
for all $x,y\in \mathfrak{g}$ and $f,g\in \mathfrak{g}^*$, where $\phi(x, f)\cdot g$ is defined by
\begin{equation}\label{action}\langle y,\phi(x, f)\cdot g\rangle=-\langle[y,\phi(x, f)], g\rangle,\ \textrm{for\ all} \ y\in\mathfrak{g}.
\end{equation}

The fundamental identity \eqref{aaa} suggests defining a bracket on $\mathfrak{g}\otimes\mathfrak{g}^*$
by
\begin{equation}\label{bbb}[x\otimes f,y\otimes g]=[\phi(x, f), y]\otimes g+y\otimes \phi(x, f)\cdot g.
\end{equation}
\begin{proposition}\cite{def3alg}The bracket \eqref{bbb} turns $\mathfrak{g}\otimes\mathfrak{g}^*$
into a  Leibniz algebra.
\end{proposition}

\begin{proposition}\label{prop}
Let $\alpha\in Aut_S(B,\mathfrak{g})$ be an involution, then $(\mathfrak{g}\otimes\mathfrak{g}^*,[.,.]_\Omega,\Omega,B_\Omega)$, where
\begin{align}
&\Omega(x\otimes f)=\alpha(x)\otimes f\circ\alpha,\\
&[x\otimes f,y\otimes g]_\Omega=\Omega\circ[x\otimes f,y\otimes g],\\
&B_\Omega(x\otimes f,y\otimes g)=\langle\alpha(x),g\rangle\langle\alpha(y),f\rangle,
\end{align}
 is a multiplicative quadratic Hom-Leibniz algebra.
\end{proposition}
\begin{proof}
Let $x,y,z\in \mathcal{N}$, $f,g,h\in\mathcal{N}^*$.  Using \eqref{sc} and \eqref{action},  we have
\begin{align*}
B(y,\alpha(\phi(x\otimes f)))&=B(\alpha(y),\phi(x\otimes f))\\
&=\langle[\alpha(y),x],f\rangle\\
&=\langle\alpha([y,\alpha(x)]),f\rangle\\
&=\langle[y,\alpha(x)],f\circ\alpha\rangle\\
&=B(y,\phi(\alpha(x)\otimes f\circ\alpha),
\end{align*}
and
\begin{align*}
\langle y,(\phi(x, f)\cdot g)\circ\alpha\rangle&=\langle \alpha(y),\phi(x, f)\cdot g\rangle\\
&=-\langle[\alpha(y),\phi(x, f)], g\rangle\\
&=-\langle\alpha([y,\alpha(\phi(x, f))]), g\rangle\\
&=-\langle[y,\alpha(\phi(x, f))], g\circ\alpha\rangle\\
&=-\langle[y,\phi(\alpha(x), f\circ\alpha)], g\circ\alpha\rangle\\
&=\langle y,\phi(\alpha(x), f\circ\alpha)\cdot (g\circ\alpha)\rangle.
\end{align*}
Thus, we obtain the following identity
\begin{align*}
&\alpha(\phi(x\otimes f))=\phi(\alpha(x)\otimes f\circ\alpha),\\
&(\phi(x, f)\cdot g)\circ\alpha=\phi(\alpha(x), f\circ\alpha)\cdot (g\circ\alpha).
\end{align*}
Therefore
\begin{align*}
\Omega([x\otimes f,y\otimes g])&=\alpha([\phi(x, f), y])\otimes g\circ\alpha+\alpha(y)\otimes (\phi(x, f)\cdot g)\circ\alpha\\
&=[\alpha(\phi(x, f)), \alpha(y)]\otimes g\circ\alpha+\alpha(y)\otimes (\phi(x, f)\cdot g)\circ\alpha\\
&=[\phi(\alpha(x)\otimes f\circ\alpha), \alpha(y)]\otimes g\circ\alpha+\alpha(y)\otimes \phi(\alpha(x), f\circ\alpha)\cdot (g\circ\alpha)\\
&=[\alpha(x)\otimes f\circ\alpha,\alpha(y)\otimes g\circ\alpha]\\
&=[\Omega(x\otimes f),\Omega(y\otimes g)].
\end{align*}
Thus, $\Omega([x\otimes f,y\otimes g])=[\Omega(x\otimes f),\Omega(y\otimes g)]$.  Then using Theorem \ref{twistcentral1}, $(\mathfrak{g}\otimes\mathfrak{g}^*,[.,.]_\Omega,\Omega)$ is a multiplicative  Hom-Leibniz algebra.

Since $\alpha$ is an involution, then
\begin{align*}
B_\Omega([x\otimes f,y\otimes g]_\Omega,z\otimes h)&=B_\Omega(\alpha([\phi(x, f), y])\otimes g\circ \alpha,z\otimes h)+B_\Omega(\alpha(y)\otimes (\phi(x, f)\cdot g)\circ\alpha,z\otimes h)\\
&=\langle [\phi(x, f), y],h\rangle\langle\alpha(z),g\circ\alpha\rangle+\langle y,h\rangle\langle\alpha(z),(\phi(x, f)\cdot g)\circ\alpha\rangle\\
&=\langle [\phi(x, f), y],h\rangle\langle z,g\rangle+\langle y,h\rangle\langle z,\phi(x, f)\cdot g\rangle\\
&=\langle [\phi(x, f), y],h\rangle\langle z,g\rangle-\langle y,h\rangle\langle[z,\phi(x, f)], g\rangle\\
&=-(\langle[z,\phi(x, f)], g\rangle\langle y,h\rangle-\langle z,g\rangle\langle [\phi(x, f), y],h\rangle)\\
&=-B_\Omega(y\otimes g,[x\otimes f,z\otimes h]_\Omega).
\end{align*}
Finally, the bilinear form $B_\Omega$ is symmetric nondegenerate and invariant.  Then  $(\mathfrak{g}\otimes\mathfrak{g}^*,[.,.]_\Omega,\Omega,B_\Omega)$ is a multiplicative quadratic  Hom-Leibniz algebra.
\end{proof}
The inner product on $\mathfrak{g}$ sets up an
isomorphism $\flat:\mathfrak{g}\rightarrow\mathfrak{g}^*$
of $\mathfrak{g}$-modules, defined by $x^*=\flat(x)=B(x,\cdot)$.\\
The map $\phi$ defined by equation \eqref{sc} induces a map $T:\mathfrak{g}\otimes\mathfrak{g}\rightarrow\mathfrak{g},$ by $T(x\otimes y)=\phi(x\otimes y^*)$. In other words, for all $x,y,z\in\mathfrak{g}$, we have
$$B(T(x\otimes y),z)=B([z,x],y),$$
whence
$$T(x\otimes y)=-T(y\otimes x ).$$
This means that $T$ factors through a map also denoted $T:\wedge^2\mathfrak{g}\rightarrow\mathfrak{g}$.\\
Using $T$ we can define a ternary bracket on $\mathfrak{g}$ by
\begin{equation}\label{Tbracket}
[x,y,z]:=[T(x\otimes y),z]
\end{equation}
and $(\mathfrak{g},[.,.,.],B)$ is a quadratic ternary Nambu algebra. 

\begin{proposition}
Let $(\mathfrak{g},[\cdot,\cdot ],B)$ be a real finite-dimensional quadratic  Lie algebra and $\alpha\in Aut_S(B,\mathfrak{g})$ be an involution. Then $(\mathfrak{g},\alpha\circ[\cdot,\cdot,\cdot],(\alpha,\alpha),B_\alpha)$,  where the bracket is defined in \ref{Tbracket} and $B_\alpha(x,y)=B(\alpha(x),y),$
 is a quadratic multiplicative ternary Hom-Nambu algebra.
 \end{proposition}

\section{Centroids, derivations and quadratic $n$-ary Hom-Nambu  algebras }
In this section we first  generalize to $n$-ary Hom-Nambu algebras the notion of centroid and its properties discussed in \cite{n-centroid}. We also generalize to Hom setting the connections between centroid elements and derivations. Finally we construct quadratic $n$-ary Hom-Nambu algebras involving elements of the centroid of $n$-ary Nambu algebras.

\subsection{Centroids of $n$-ary Hom-Nambu  algebras }
\begin{definition}
Let  $(\mathcal{N},[\cdot,...,\cdot],\alpha)$ be a multiplicative $n$-ary Hom-Nambu algebra and $End(\mathcal{N})$ be the endomorphism  algebra of $\mathcal{N}$. Then the following subalgebra
of $End(\mathcal{N})$
\begin{equation}
Cent(\mathcal{N})=\{\theta\in End(\mathcal{N}): \theta[x_1,...,x_n] = [\theta x_1,...,x_n],\ \forall x_i\in\mathcal{N}\}.
\end{equation}
is said to be the  centroid of the $n$-ary Hom-Nambu algebra.

The definition is the same for  classical case of $n$-ary Nambu algebra. We may also consider the same definition for any $n$-ary Hom-Nambu algebra.
\end{definition}

Now, let $(\mathcal{N},[\cdot,...,\cdot],\alpha)$ be a multiplicative $n$-ary Hom-Nambu algebra. We denote by $\alpha^k$, where $\alpha \in End(\mathcal{N})$,  the $k$-times composition of $\alpha$. We set in
 particular $\alpha^{-1}=0$ and $\alpha^0=Id$.

\begin{definition}

An  $\alpha^k$-centroid of a multiplicative $n$-ary Hom-Nambu algebra $(\mathcal{N},[\cdot,...,\cdot],\alpha)$ is a subalgebra
of $End(\mathcal{N})$ denoted $Cent_{\alpha^k}(\mathcal{N})$, given by
\begin{equation}
Cent_{\alpha^k}(\mathcal{N})=\{\theta\in End(\mathcal{N}):\theta[x_1,...,x_n] = [\theta x_1,\alpha^k(x_2)...,\alpha^k(x_n)],\ \forall x_i\in\mathcal{N}.
\end{equation}
\end{definition}
We recover the definition of the centroid when $k=0$.

If $\mathcal{N}$ is a multiplicative $n$-ary Hom-Nambu-Lie algebra, then it is a simple fact that
$$\theta[x_1,....,x_n]=[\alpha^k(x_1),...,\theta x_p,...,\alpha^k(x_n)],\ \ \forall \ p\in\{1,...,n\}.$$

\begin{lem}\label{centalg}
Let $(\mathcal{N}, [\cdot,...,\cdot ])$ be an $n$-ary Nambu-Lie algebra. If $\theta\in Cent(\mathcal{N})$, then for $x_1,...,x_n\in \mathcal{N}$

\begin{enumerate}
  \item $[\theta^{p_1} x_1,...,\theta^{p_n} x_n]=\theta^{p_1+...+p_n}[ x_{1},... ,x_{ n}], \ \forall p_1,...,p_n\in
  \mathbb{N}$,
  \item $[\theta^{p_1} x_1,...,\theta^{p_n} x_n]=Sgn(\sigma) [\theta^{p_1} x_{\sigma(1)},...,\theta^{p_n} x_{\sigma (n)}],  \ \forall p_1,...,p_n\in \mathbb{N}$
 and $\forall  \sigma\in\mathcal{S}_n$.

\end{enumerate}
\end{lem}

\begin{proof}
Let $\theta\in Cent(\mathcal{N})$, $x_1,...,x_n\in \mathcal{N}$ and $1\leq p\leq n$, we have
$$[\theta^p x_1,...,x_n]=\theta[\theta^{p-1}x_1,...,x_n]=...=\theta^p[x_1,...,x_n].$$
Also, observe that for any $k\in\{1,...,n\}$
$$ [x_1,...,\theta^p x_k,...,x_n]=-[\theta^p x_k,x_2,...,x_1,...,x_n]=-\theta^p[x_k,x_2,...,x_1,...,x_n]=\theta^p[x_1,...,x_k,...,x_n].$$
Then, similarly we have
\begin{align*}
&[\theta^{p_1} x_1,...,\theta^{p_n} x_n]=\theta^{p_n}[\theta^{p_1} x_1,...,\theta^{p_{n-1}} x_{n-1}, x_n]=...=\theta^{p_1+...+p_n}[ x_1,...,x_n].
\end{align*}
The second assertion is a consequence of previous calculations and the skew-symmetry of $[\cdot,...,\cdot]$.
\end{proof}

\begin{proposition}\label{centconstruction}
Let $(\mathcal{N}, [\cdot,...,\cdot ])$ be an $n$-ary Nambu-Lie algebra and $\theta\in Cent(\mathcal{N})$.\\ Let us fix $p$ and set for any $x_1,...,x_n \in\mathcal{N}$
\begin{equation}\{x_1,...,x_n\}_ p= [\theta x_1,...,\theta x_{p-1},\theta x_p,x_{p+1},...,x_n].\end{equation}
Then $(\mathcal{N}, \{\cdot,...,\cdot \}_p,\widetilde{\theta}=(\theta,...,\theta))$ is an $n$-ary Hom-Nambu-Lie algebra.
\end{proposition}

\begin{proof}
For $\theta\in Cent(\mathcal{N})$ and $p\in\{1,...,n\}$,  we have
\begin{align*}
\{\theta x_1,...,\theta x_{n-1},\{y_1,...,y_n\}_p\}_p&= [\theta^2 x_1,...,\theta^2 x_p,...,\theta x_{n-1},[\theta y_1,...,\theta y_p,...,y_n]]\\
&=[\theta^2 x_1,...,\theta^2 x_p,...,\theta x_{n-1},\theta^p[ y_1,...,y_n]]\\
&=\theta^{2p+n-1}([ x_1,..., x_{n-1},[ y_1,...,y_n]]).
\end{align*}
In the other hand we have
\begin{align*}
&\sum_{k=0}^n\{\theta y_1,...,\{ x_1,..., x_{n-1},y_k\}_p,...,\theta y_n\}_p\\
&= \sum_{k=0}^p\{\theta y_1,...,\{ x_1,..., x_{n-1},y_k\}_p,...,\theta y_n\}_p+\sum_{k=p}^n\{\theta y_1,...,\{ x_1,..., x_{n-1},y_k\}_p,...,\theta y_n\}_p\\
&= \sum_{k=0}^p[\theta^2 y_1,...,\theta[ \theta x_1,...,\theta x_p,..., x_{n-1},y_k],...,\theta^2 y_p,...,\theta y_n]\\
&+\sum_{k=0}^p[\theta^2 y_1,...,\theta^2 y_p,...,[ \theta x_1,...,\theta x_p,..., x_{n-1},y_k],...,\theta y_n]\\
&= \sum_{k=0}^p\theta^{2p+n-1}[ y_1,...,[  x_1,..., x_{n-1},y_k],..., y_n]+\sum_{k=0}^p\theta^{2p+n-1}[2 y_1,...,[  x_1,..., x_{n-1},y_k],..., y_n]\\
&= \theta^{2p+n-1}(\sum_{k=0}^n[ y_1,...,[  x_1,..., x_{n-1},y_k],..., y_n]).
\end{align*}
Therefore
the Hom-Nambu identity with respect to the bracket $[\cdot,...,\cdot ]$ leads to the Hom-identity for  $\{\cdot,...,\cdot \}_l$.
The skew-symmetry
is proved by second assertion of Lemma \ref{centalg}.
\end{proof}

\subsection{Centroids and Derivations of Hom-Nambu algebras}
Let $(\mathcal{N}, [\cdot  ,..., \cdot ],  \alpha )$ be a multiplicative $n$-ary Hom-Nambu-Lie algebra.

\begin{definition}
For any $k\geq1$, we call $D\in End(\mathcal{N})$ an $\alpha^k$-\emph{derivation } of the
multiplicative $n$-ary Hom-Nambu-Lie  $(\mathcal{N}, [\cdot ,...,\cdot],  \alpha )$ if $D$ and $\alpha$ commute
and we have
\begin{equation}\label{alphaKderiv2}
D[x_1,...,x_n]=\sum_{i=1}^n[\alpha^k(x_1),...,\alpha^k(x_{i-1}),D(x_i),\alpha^k(x_{i+1}),...,\alpha^k(x_n)].
\end{equation}
We denote by $Der_{\alpha^k}(\mathcal{N})$ the set of $\alpha^k$-derivations.
\end{definition}

For $x=(x_1,...,x_{n-1})\in \mathcal{N}^{ n-1}$ satisfying $\alpha(x)=x$ and $k\geq 1$,
we define the map $ad_k(x)\in End(\mathcal{N})$ by
\begin{equation}\label{ad_k(u)}
ad_k(x)(y)=[x_1,...,x_{n-1},\alpha^k(y)]\ \ \forall y\in \mathcal{N}.
\end{equation}

The map $ad_k(x)$ is an $\alpha^{k+1}$-derivation, that we call inner $\alpha^{k+1}$-derivation.
We denote by $Inn_{\alpha^k}(\mathcal{N})$ the space generated by all  inner $\alpha^{k+1}$-derivations.

Set $Der(\mathcal{N})=\dl\bigoplus_{k\geq -1}Der_{\alpha^k}(\mathcal{N})$ and $Inn(\mathcal{N})=\dl\bigoplus_{k\geq -1}Inn_{\alpha^k}(\mathcal{N})$.
\begin{lem}\label{2.3}
For $D\in Der_{\alpha^k}(\mathcal{N})$ and $D'\in Der_{\alpha^{k'}}(\mathcal{N})$, where $k+k'\geq-1$, we have
$[D,D']\in Der_{\alpha^{k+k'}}(\mathcal{N})$, where the commutator $[D,D']$ is defined as usual.
\end{lem}
Now, we define a linear map $\varsigma:Der_{\alpha^k}(\mathcal{N})\rightarrow Der_{\alpha^{k+1}}(\mathcal{N})$ by $\varsigma(D)=\alpha\circ D.$
Since  the elements of $Der_{\alpha^k}(\mathcal{N}) $ and $\alpha$ commute then $\varsigma$ is in the  centroid  of the Lie algebra
$(Der(\mathcal{N}),[\cdot ,\cdot ])$.\\Hence, using  Proposition \ref{centconstruction} we have
\begin{proposition}Let $(\mathcal{N}, [\cdot  ,..., \cdot ],  \alpha )$ be a multiplicative $n$-ary Hom-Nambu-Lie algebra.
The triple $(Der(\mathcal{N}),[\cdot ,\cdot ]_\varsigma,\varsigma)$, where the bracket is defined by $[\cdot ,\cdot ]_\varsigma=\varsigma\circ [\cdot ,\cdot ]$,  is a Hom-Lie algebra.
\end{proposition}

\begin{proposition}\label{centroidDer}Let $(\mathcal{N}, [\cdot  ,..., \cdot ],  \alpha )$ be a multiplicative $n$-ary Hom-Nambu-Lie algebra.
If $D\in Der_{\alpha^k}(\mathcal{N})$ and $\theta\in Cent_{\alpha^{k'}}(\mathcal{N})$, then $\theta D\in Der_{\alpha^{k+k'}}(\mathcal{N})$.
\end{proposition}

\begin{proof}
Let $x_1,...,x_n\in \mathcal{N}$ then
\begin{align*}
\theta D([x_1,...,x_n])&=\sum_{i=1}^n\theta [\alpha^k(x_1),...,D(x_i),...,\alpha^k(x_n)]\\
&=\sum_{i=1}^n [ \alpha^{k+k'}(x_1),...,\theta D(x_i),...,\alpha^{k+k'}(x_n)].
\end{align*}
Thus $\theta D$ is an $\alpha^k$-derivation.
\end{proof}

Now we define the notion of central derivation. Let $(\mathcal{N}, [\cdot  ,..., \cdot ],  \alpha )$ be a multiplicative $n$-ary Hom-Nambu-Lie algebra. We set  $Z(\mathcal{N})=\{x\in \mathcal{N}:[x,y_1,...,y_{n-1}]=0,\ \forall \ y_1,...,y_{n-1}\in\mathcal{N}\}$, the center of the $n$-ary Hom-Nambu-Lie algebra.
\begin{definition}
Let  $\varphi\in End(\mathcal{N})$, then  $\varphi$ is said to be a central
derivation if $\varphi(\mathcal{N}) \subset Z(\mathcal{N})$ and $\varphi([\mathcal{N},...,\mathcal{N}]) = 0$.

The set of all central derivations of $\mathcal{N}$ is denoted by $C(\mathcal{N})$.
\end{definition}
Notice that  an  $\alpha^k$-derivation $\varphi$ is a central derivation if $\varphi(\mathcal{N}) \subset Z(\mathcal{N})$.
\begin{thm}
Let $(\mathcal{N}, [\cdot  ,..., \cdot ],  \alpha )$ be a multiplicative $n$-ary Hom-Nambu-Lie algebra.
Let  $D$ in $Der_{\alpha^k}(\mathcal{N})$ and $\theta$ in $Cent(\mathcal{N})$ such that $[\theta,\alpha]=0$, then we have
\begin{enumerate}
  \item $[D,\theta]$ is in the  $\alpha^k$-centroid of $\mathcal{N}$,
  \item if $ [D,\theta]$ is a central derivation then $D\theta$ is an $\alpha^k$-derivation of $\mathcal{N}$.
\end{enumerate}
\end{thm}

\begin{proof}(1)
Let $D\in Der_{\alpha^k}(\mathcal{N})$, $\theta\in Cent(\mathcal{N})$ and $x_1,...,x_n\in \mathcal{N}$ we have
\begin{align*}
D\theta([x_1,...,x_n])&=D([\theta x_1,...,x_n])\\
&=[D\theta x_1,...,\alpha^k(x_n)]+ \sum_{i=2}^n [\alpha^k(\theta x_1),...,D(x_i),...,\alpha^k(x_n)]\\
&=[D\theta x_1,...,\alpha^k(x_n)]+ \sum_{i=2}^n [\theta\alpha^k( x_1),...,D(x_i),...,\alpha^k(x_n)]\\
&=[D\theta x_1,...,\alpha^k(x_n)]+ \sum_{i=2}^n [\alpha^k( x_1),...,\theta D(x_i),...,\alpha^k(x_n)]\\
&=[D\theta x_1,...,\alpha^k(x_n)]+ \theta D([x_1,...,x_n])-[\theta D x_1,...,\alpha^k(x_n)].
\end{align*}
Then
$$(D\theta- \theta D)([x_1,...,x_n])=[(D\theta -\theta D) x_1,\alpha^k(x_2),...,\alpha^k(x_n)].$$
That is, $[D,\theta]=D\theta-\theta D\in Cent_k(\mathcal{N})$.

(2) Using Proposition \ref{centroidDer}, $\theta D$ is an $\alpha^k$-derivation and since $[D,\theta] $ is a $\alpha^k$-derivation, then $D\theta=[D,\theta]+\theta D$ is also an $\alpha^k$-derivation.
\end{proof}

Let $A$ be a $\K$-vector space, $\mu$ be an $n$-linear map on $A$ and $\eta$  be  linear maps on $A$. Let $(A,\mu,\eta)$  be a multiplicative symmetric $n$-ary totally Hom-associative algebra. The $\eta^k$-centroid $Cent_{\eta^k}(A)$ of $A$ is  defined by
$$Cent_{\eta^k}(A)=\{f\in End(A): \ f(\mu(a_1,...,a_n))=\mu(f(a_1),\eta^k(a_2),..., \eta^k(a_n))\},$$for all $ a_i\in A \ \textrm{and}\ i\in\{1,...,n\}$.
The set of $\eta^k$-derivation, $Der_{\eta^k}(A)$, is a subset of  $End(A)$ defined by $\varphi\in End(A)$
such that $\varphi(\mu(a_1,...,a_n))=\dl\sum_{i=1}^n\mu(\eta^k(a_1)...., \eta^k(a_{i-1}),\varphi(a_i), \eta^k(a_{i+1}),..., \eta^k(a_n))\},$ for all $ a_i\in A$.

\begin{thm}Let $(A,\mu,\eta)$  be a  multiplicative symmetric $n$-ary Hom-associative algebra and $(\mathcal{N}, [\cdot  ,..., \cdot ],  \alpha )$ be a  multiplicative $n$-ary Hom-Nambu-Lie algebra, then we have the following assertion
\begin{itemize}

      \item If $f\in Cent_{\eta^k}(A)$ and $\theta\in Cent_{\alpha^k}(\mathcal{N})$, then $f\otimes \theta  $ is in the  $\zeta^k$-centroid, where $\zeta^k=\eta^k\otimes\alpha^k$,  of the Hom-Nambu-Lie algebra $A\otimes\mathcal{N}$ defined in  \ref{prodtens}.
          \item If $f\in Cent_{\eta^k}(A)$ and $D\in Der_{\alpha^k}(\mathcal{N})$, then $f\otimes D  $ is a $\zeta^k$-derivation of the Hom-Nambu-Lie algebra $A\otimes\mathcal{N}$.
\end{itemize}
\end{thm}

\begin{proof}Let $a_i\in A$, $x_i\in \mathcal{N}$ where $i\in\{1,...,n\}$ and $f$ be a $\eta^k$-centroid on $A$. \\
$\bullet$ If $\theta\in Cent_{\alpha^k}(\mathcal{N})$, then
\begin{align*}
(f\otimes \theta)( [a_1\otimes x_1,...,a_n\otimes x_n])&=(f\otimes \theta) (\mu(a_1,..., a_n)\otimes[ x_1,..., x_n]_\mathcal{N})\\
&=  \mu(f(a_1),\eta^k(a_2),...,\eta^k( a_n))\otimes[ \theta(x_1),\alpha^k(x_2),...,\alpha^k( x_n)]_\mathcal{N}\\
&= [(f\otimes \theta)(a_1\otimes x_1),\zeta^k(a_2\otimes x_2),...,\zeta^k(a_n\otimes x_n)].
\end{align*}
Thus $f\otimes \theta$ is in the  $\zeta^k$-centroid of $A\otimes\mathcal{N}$.\\
$\bullet$ If $D\in Der_{\alpha^k}(\mathcal{N})$, then
\begin{align*}
(f\otimes D)( [a_1\otimes x_1,...,a_n\otimes x_n])&=f\otimes D((a_1\cdot...\cdot a_n)\otimes[ x_1,..., x_n])\\
&=  \mu(f(a_1),\eta^k(a_2),...,\eta^k( a_n))\otimes\sum_{i=1}^n[ \alpha^k(x_1),...,D(x_i),...,\alpha^k( x_n)]\\
&= \sum_{i=1}^n\mu(\eta^k(a_1),..., f(a_i),...,\eta^k( a_n))\otimes[ \alpha^k(x_1),...,D(x_i),...,\alpha^k( x_n)]_\mathcal{N}\\
&=\sum_{i=1}^n[\zeta^k(a_1\otimes x_1),...,(f\otimes D)(a_i\otimes x_i),...,\zeta^k(a_n\otimes x_n)].
\end{align*}
Therefore $f\otimes D$ is a $\zeta^k$-derivation of $A\otimes\mathcal{N}$.
\end{proof}

\subsection{ Centroids and Quadratic $n$-ary Hom-Nambu algebras}
Let $\theta\in Cent(\mathcal{N})$ such that $\theta$ is invertible and symmetric with respect to  $B$ $($i.e.  $B(\theta x,y)=B(x,\theta y)$ $)$.
We set
$$Cent_S(\mathcal{N})=\{\theta\in Cent(\mathcal{N}): \theta \ \textrm{symmetric with respect to }\ B\}.$$
\begin{thm}
Let $(\mathcal{N},[\cdot,...,\cdot],B)$
be a quadratic $n$-ary Nambu-Lie algebra and $\theta\in Cent_S(\mathcal{N})$ such that $\theta$ is invertible. We consider  a  bilinear form  $B_\theta$ defined  by
$$\begin{array}{cccc}
  B_\theta: & \mathcal{N}\times \mathcal{N}&\longrightarrow& \K \\
    & (x,y) &\longmapsto & B(\theta x,y) .
\end{array}$$
Then, $(\mathcal{N},\{\cdot,...,\cdot\}_l,(\theta,...,\theta),B_\theta)$ is a quadratic $n$-ary Hom-Nambu-Lie algebra.
\end{thm}
\begin{proof}
It easy to proof that $B_\theta$ is symmetric and nondegenerate.\\We have also $\theta$ is symmetric with respect to $B_\theta$, indeed
$$B_\theta(\theta x,y)=B(\theta^2x,y)=B(\theta x,\theta y)=B_\theta( x,\theta y).$$
The invariance of $B_\theta$ is given by, set $l\in\{1,...,n-1\}$
\begin{align*}
B_\theta(\{x_1,...,x_{n-1},y\}_l,z)&=B_\theta([\theta x_1,...,\theta x_l,...,x_{n-1},y],z)\\
&=B(\theta[\theta x_1,...,\theta x_l,...,x_{n-1},y],z)\\
&=B([\theta^2 x_1,...,\theta x_l,...,x_{n-1},y],z)\\
&=-B(y,[\theta^2 x_1,...,\theta x_l,...,x_{n-1},z])\\
&=-B_\theta(y,[\theta x_1,...,\theta x_l,...,x_{n-1},z])\\
&=-B_\theta(y,\{x_1,...,x_{n-1},z\}_l)
\end{align*}
In the other hand, when  $l=n$ we have $[\theta x_1,...,\theta x_n]=[\theta^2 x_1,...,\theta x_{n-1},  x_n]$ then it's a consequence of the previous calculations.\\
Therefore  $(\mathcal{N},\{\cdot,...,\cdot\}_l,(\theta,...,\theta),B_\theta)$ is a quadratic $n$-ary Hom-Nambu-Lie algebra.\\
Notice that $B_\theta $ is also an invariant scalar product of the $n$-ary Nambu-Lie algebra $\mathcal{N}$.
\end{proof}

\end{document}